\documentclass{article}%
\usepackage{graphicx}
\usepackage{amsmath}
\usepackage{amsfonts}
\usepackage{amssymb}%
\newtheorem{theorem}{Theorem}[section]

\newtheorem{lemma}[theorem]{Lemma}
\newtheorem{proposition}[theorem]{Proposition}
\newtheorem{definition}[theorem]{Definition}

\newenvironment{proof}[1][Proof]{\textbf{#1.} }{\ \rule{0.5em}{0.5em}}
\begin{document}

\title{On finitely generated profinite groups, II: products in quasisimple groups}
\author{Nikolay Nikolov\thanks{Work done while the first author held a Golda-Meir
Fellowship at the Hebrew University of Jerusalem} \ and Dan Segal}
\maketitle

\begin{abstract}
We prove two results. (1) There is an absolute constant $D$ such that for any
finite quasisimple group $S$, given $2D$ arbitrary automorphisms of $S$, every
element of $S$ is equal to a product of $D$ `twisted commutators' defined by
the given automorphisms.

(2) Given a natural number $q$, there exist $C=C(q)$ and $M=M(q)$ such that:
if $S$ is a finite quasisimple group with $\left|  S/\mathrm{Z}(S)\right|
>C$, $\beta_{j}$ $\ (j=1,\ldots,M)$ are any automorphisms of $S$, and $q_{j}$
$\ (j=1,\ldots,M)$ are any divisors of $q$, then there exist inner
automorphisms $\alpha_{j}$ of $S$ such that $S=\prod_{1}^{M}[S,(\alpha
_{j}\beta_{j})^{q_{j}}]$.

These results, which rely on the Classification of finite simple groups, are
needed to complete the proofs of the main theorems of Part I.

\end{abstract}

\section{Introduction}

The main theorems of Part I \cite{NS} were reduced to two results about finite
quasisimple groups. These will be proved here.

A group $S$ is said to be \emph{quasisimple} if $S=[S,S]$ and $S/\mathrm{Z}%
(S)$ is simple, where $\mathrm{Z}(S)$ denotes the centre of $S$. For
automorphisms $\alpha$ and $\beta$ of $S$ we write
\[
T_{\alpha,\beta}(x,y)=x^{-1}y^{-1}x^{\alpha}y^{\beta}.
\]

\begin{theorem}
\label{twisted}There is an absolute constant $D\in\mathbb{N}$ such that if $S$
is a finite quasisimple group and $\alpha_{1},\,\beta_{1},\ldots,\alpha
_{D},\,\beta_{D}$ \ are any automorphisms of $S$ then
\[
S=T_{\alpha_{1},\beta_{1}}(S,S)\cdot\cdots\cdot T_{\alpha_{D},\beta_{D}%
}(S,S).
\]

\end{theorem}

\begin{theorem}
\label{autos}Let $q$ be a natural number. There exist natural numbers $C=C(q)$
and $M=M(q)$ such that if $S$ is a finite quasisimple group with $\left|
S/\mathrm{Z}(S)\right|  >C$, $\beta_{1},\ldots,\beta_{M}$ \ are any
automorphisms of $S$, and $q_{1},\ldots,q_{M}$ are any divisors of $q$, then
there exist inner automorphisms $\alpha_{1},\ldots,\alpha_{M}$ of $S$ such
that
\[
S=[S,(\alpha_{1}\beta_{1})^{q_{1}}]\cdot\cdots\cdot\lbrack S,(\alpha_{M}%
\beta_{M})^{q_{M}}].
\]

\end{theorem}

These results are stated as Theorems 1.9 and 1.10 in the introduction of
\cite{NS}. Both may be seen as generalizations of Wilson's theorem \cite{W}
that every element of any finite simple group is equal to the product of a
bounded number of commutators: indeed, we shall show that in Theorem
\ref{autos}, $C(1)$ may be taken equal to $1$. The latter theorem also
generalizes the theorem of Martinez, Zelmanov, Saxl and Wilson (\cite{MZ},
\cite{SW}) that in any finite simple group $S$ with $S^{q}\neq1$, every
element is equal to a product of boundedly many $q$th powers, the bound
depending only on $q$.

The proofs depend very much on the Classification of finite simple groups, and
in Section 2 we give a brief r\'{e}sum\'{e} on groups of Lie type, its main
purpose being to fix a standard notation for these groups, their subgroups and
automorphisms. Section 3 collects some combinatorial results that will be used
throughout the proof, and in Section 4 we show that Theorem \ref{twisted} is a
corollary of Theorem \ref{autos} (it only needs the special case where $q=1$).

The rest of the paper is devoted to the proof of Theorem \ref{autos}. This
falls into two parts. The first, given in \S 5, concerns the case where
$S/\mathrm{Z}(S)$ is either an alternating group or a group of Lie type over a
`small' field; this case is deduced from known results by combinatorial
arguments. The second part, in \S \S 6--10, deals with groups of Lie type over
`large' fields: this depends on a detailed examination of the action of
automorphisms of $S$ on the root subgroups. The theorem follows, since
according to the Classification all but finitely many of the finite simple
groups are either alternating or of Lie type.

\subsection*{Notation}

For a group $S$ and $x,\,y\in S$, $x^{y}=y^{-1}xy$. If $y\in S$ or
$y\in\mathrm{Aut}(S)$,%
\[
\lbrack x,y]=x^{-1}x^{y},\;[S,y]=\left\{  [x,y]\mid x\in S\right\}  .
\]

We will write $\overline{S}=S/\mathrm{Z}(S)$, and identify this with the group
$\mathrm{Inn}(S)$ of inner automorphisms of $S$. Similarly for $g \in S$ we
denote by $\bar{g} \in\mathrm{Inn}(S)$ the automorphism induced by conjugation
by $g$. The Schur multiplier of $S$ is denoted $M(S)$, and $\mathrm{Out}(S)$
denotes the outer automorphism group of $S$. For a subset $X\subseteq S$, the
subgroup generated by $X$ is denoted $\left\langle X\right\rangle ,$ and for
$n\in\mathbb{N}$%
\[
X^{\ast n} =\left\{  x_{1}\ldots x_{n}\mid x_{1},\ldots,x_{n}\in X\right\}  .
\]

We use the usual notation $F^{\ast}$ for the multiplicative group
$F\setminus\{0\}$ of a field $F$ (this should cause no confusion). The symbol
$\log$ means logarithm to base $2$.

\section{Groups of Lie type: a r\'{e}sum\'{e}\label{resume}}

Apart from the alternating groups $\mathrm{Alt}(k)$ $(k\geq5)$ and finitely
many sporadic groups, every finite simple group is a \emph{group of Lie type},
that is, an untwisted or twisted Chevalley group over a finite field. We
briefly recall some features of these groups, and fix some notation. Suitable
references are Carter's book \cite{C},  Steinberg's lectures \cite{St} and
\cite{GLS}. For a useful summary (without proofs) see also Chapter 3 of
\cite{At}.

\subsection*{\medskip Untwisted Chevalley groups}

Let $\mathcal{X}$ be one of the Dynkin diagrams $A_{r}\ (r\geq1),B_{r}%
\ (r\geq2),C_{r}\ (r\geq3),D_{r}\ (r\geq4),E_{6},E_{7},E_{8},F_{4}$ or $G_{2}%
$, let $\Sigma$ be an irreducible root system of type $\mathcal{X}$ and let
$\Pi$ be a fixed base of fundamental roots for $\Sigma$. This determines
$\Sigma_{+}$: the positive roots of $\Sigma$.

The number $r:=\left|  \Pi\right|  $ of fundamental roots is the \emph{Lie
rank}. If $w\in\Sigma_{+}$ is a positive root then we can write $w$ uniquely
as a sum of fundamental roots (maybe with repetitions). The number of
summands, denoted $\mathrm{ht}(w)$, is called the \emph{height} of $w$. Thus
$\Pi$ is exactly the set of roots of height $1$.

Let $\mathcal{X}(F)$ denote the $F$-rational points of a split simple
algebraic group of type $\mathcal{X}$ over the field $F$. To each $w\in\Sigma$
is associated a one-parameter subgroup of $\mathcal{X}(F)$,
\[
X_{w}=\left\{  X_{w}(t)\mid t\in F\right\}  ,
\]
called the \emph{root subgroup} corresponding to $w$.

The associated \emph{Chevalley group} of type $\mathcal{X}$ over $F$ is
defined to be the subgroup $S$ of $\mathcal{X}(F)$ generated by all the root
subgroups $X_{w}$ for $w\in\Sigma$. It is \emph{adjoint }%
(resp\emph{\ universal}) if the algebraic group is adjoint (resp. simply
connected). With finitely many exceptions $S$ is a quasisimple group.

Let $U=U_{+}:=\prod_{w\in\Sigma_{+}}X_{w}$ and $U_{-}:=\prod_{w\in\Sigma_{-}%
}X_{w}$, the products being ordered so that $\mathrm{ht}(\left|  w\right|  )$
is non-decreasing. Then $U_{+}$ and $U_{-}$ are subgroups of $S$ (the
positive, negative, unipotent subgroups).

For each multiplicative character $\chi:\mathbb{Z}\Sigma\rightarrow F^{\ast}$
of the lattice spanned by $\Sigma$, we define an automorphism $h(\chi)$ of $S$
by
\[
X_{w}(t)^{h(\chi)}=X_{w}(\chi(w)\cdot t).
\]
The set of all such $h(\chi)$ forms a subgroup $\mathcal{D}$ of $\mathrm{Aut}%
(S)$, called the group of \emph{diagonal automorphisms}.

The group $H$ of \emph{diagonal elements} of $S$ is the subgroup generated by
certain semisimple elements $h_{v}(\lambda)$ ($v\in\Pi,\,\lambda\in F^{\ast}%
$); the group $H$ normalizes each root subgroup, and we have
\[
X_{w}(t)^{h_{v}(\lambda)}=X_{w}(\lambda^{\langle w,v\rangle}t)
\]
where $\langle w,v\rangle=2(w,v)/\left|  v\right|  ^{2}$. In particular
$X_{w}(t)^{h_{w}(\lambda)}=X_{w}(\lambda^{2}t)$. The inner automorphisms of
$S$ induced by $H$ are precisely the inner automorphisms lying in
$\mathcal{D}$, and $\mathcal{D}$ acts trivially on $H$.

For each power $p^{f}$ of $\mathrm{char}(F)$ there is a \emph{field
automorphisms} $\phi=\phi(p^{f})$ of $S$ defined by
\[
X_{w}(t)^{\phi}=X_{w}(t^{p^{f}}).
\]
The set $\Phi$ of field automorphisms is a group isomorphic to $\mathrm{Aut}%
(F)$.

The groups $\mathcal{D}$ and $\Phi$ stabilize each root subgroup and each of
the \emph{diagonal subgroups} $H_{v}=\{h_{v}(\lambda)\mid\lambda\in F^{\ast
}\}$.

We write $\mathrm{Sym}(\mathcal{X})$ for the group of (root-length preserving)
symmetries of the root system $\Sigma$. This is a group of order at most $2$
except for $\mathcal{X}=D_{4}$, in which case $\mathrm{Sym}(\mathcal{X}%
)\cong\mathrm{Sym}(3)$. Let $\tau\in\mathrm{Sym}(\mathcal{X})$ be a symmetry
that preserves the weight lattice of the algebraic group $\mathcal{X}(-)$
(e.g. if the isogeny type of $\mathcal{X}(-)$ is simply connected or adjoint).
Then (cf. Theorem 1.15.2(a) of \cite{GLS}) there exists an automorphism of
$S$, \ denoted by the same symbol $\tau$, which permutes the root subgroups in
the same way as $\tau$ acts on $\Sigma$, in fact:
\[
(X_{w}(t))^{\tau}=X_{w^{\tau}}(\epsilon_{w}t),\quad\epsilon_{w}\in
\{\pm1\}\text{ with }\epsilon_{w}=1\text{ if }w\in\Pi.
\]
This is called an \emph{ordinary graph automorphism}.

In case $\mathcal{X}=B_{2},G_{2},F_{4}$ and $p=2,3,2$ respectively, such an
automorphism of $S$ exists also when $\tau$ corresponds to the (obvious)
symmetry of order $2$ of the Dynkin diagram, which does \emph{not} preserve
root lengths. It is defined by
\[
X_{w}(t)^{\tau_{0}}=X_{w^{\tau}}(t^{r}),\quad r=1\text{ if }w\text{ is long,
}\ r=p\text{ if }w\text{ is short }.
\]
In this case $\tau_{0}$ is called an \emph{extraordinary} \emph{graph
automorphism}, and we set $\Gamma=\{1,\tau_{0}\}$. In all other cases, we
define $\Gamma$ to be the set of all ordinary graph automorphisms.

Observe that $\Gamma\neq\{1\}$ only when the rank is small ($\leq6$) or when
$\mathcal{X}$ is $A_{r}$ or $D_{r}$. The set $\Gamma$ is a group unless
$\mathcal{X}$ is one of$\ B_{2}(2^{n}),G_{2}(3^{n})$ and $F_{4}(2^{n})$, when
$|\Gamma|=2$ and the extraordinary element of $\Gamma$ squares to the
generating field automorphism of $\Phi$. In all cases, $\Gamma$ is a set of
coset representatives for $\mathrm{Inn}(S)\mathcal{D}\Phi$ in $\mathrm{Aut}%
(S)$.

\subsection*{Twisted Chevalley groups}

These are of types $^{2}A_{r},^{2}B_{2},^{2}D_{r},^{2}E_{6},^{2}F_{4}%
,^{2}G_{2},$ and $^{3}D_{4}$.

The \emph{twisted group} $S^{\ast}$ of type $^{n}\mathcal{X}$ is associated to
a certain graph automorphism of an untwisted Chevalley group $S$ of type
$\mathcal{X}$. The structure of $S^{\ast}$ is related to that of $S$, the most
notable difference being that root subgroups may be no longer one-parameter.
Another difference is that $S^{\ast}$ does not have graph automorphisms.

Let $\tau\in\Gamma\backslash\{1\}$ be a graph automorphism of $S$ as defined
above, and let $n\in\{2,3\}$ be the order of the symmetry $\tau$ on $\Sigma$.
The group $S^{\ast}$ is the fixed-point set in $S$ of the so-called
\emph{Steinberg automorphism} $\sigma=\phi\tau$, where $\phi$ is the
nontrivial field automorphism chosen so that $\sigma$ has order $n$.

We define $F_{0}\subseteq F$ to be the fixed field of $\phi$ if $\mathcal{X\,}%
$has roots of only one length; otherwise set $F_{0}=F$. In all cases,
$(F:F_{0})\leq3$.

The \emph{(untwisted) rank} of $S^{\ast}$ is defined to be the Lie rank $r$ of
the (original) root system $\Sigma$.

The root subgroups of $S^{\ast}$ now correspond to equivalence classes
$\omega$ under the equivalence relation on $\Sigma$ defined as follows.

Let $\Sigma$ be realized as a set of roots in some Euclidean vector space $V$.
The symmetry $\tau$ extends to a linear orthogonal map of $V$ and by $v^{\ast
}$ we denote the orthogonal projection of $v \in\Sigma$ on $C_{V}(\tau)$, the
subspace fixed by $\tau$. Now, for $u,v \in\Sigma$ define
\[
u \sim v \ \text{ if } \ u^{\ast}= q v^{\ast}\ \text{ for some positive } \ q
\in\mathbb{Q} .
\]
Each equivalence class $\omega$ of $\Sigma/ \sim$ is the positive integral
span of a certain orbit $\overline{\omega}$ of $\tau$ on the root system
$\Sigma$.

The root subgroups are the fixed points of $\sigma$ acting on
\[
W_{\omega}: =\langle X_{v} | \quad v \in\omega\rangle\leq S.
\]

In order to distinguish them from the root subgroups of the corresponding
untwisted group we denote them by $Y_{\omega}$. For later use we list their
structure and multiplication rules below (cf. Table 2.4 of \cite{GLS}).

Let $p=\mathrm{char}(F)$ and suppose that $\phi(t)=t^{p^{f}}$. The
automorphism $t \mapsto t^{p}$ of $F$ is denoted by $[p]$. \medskip

\textbf{Case $A_{1}^{d}$:} $\omega$ is one of $\{v\},~\{v,v^{\tau
}\},~\{v,v^{\tau},v^{\tau^{2}}\}$ and consists of pairwise orthogonal roots;
here $\left\vert \omega\right\vert =d$. When $d=2$ and there are two root
lenghts, $v$ is a long root. Then
\[
Y_{\omega}=\{Y_{\omega}(t):=\prod_{i=0}^{d-1}X_{v^{\tau^{i}}}(t^{\phi^{i}%
})\mid t\in\widetilde{F}\}
\]
is a one-parameter group; here $\widetilde{F}=F$ except when $d=1$ when
$\widetilde{F}=F_{0}$.\medskip

\textbf{Case $A_{2}$:} $\omega=\{w,v,w+v\}$ is of type $A_{2}$ with symmetry
$\tau$ swapping $v$ and $w$. Here $\left|  F\right|  =p^{2f}$, $\left|
F_{0}\right|  =p^{f},$ $\phi^{2}=1$. Elements of $Y_{\omega}$ take the form
\[
Y_{\omega}(t,u)=X_{v}(t)X_{w}(t^{\phi})X_{v+w}(u)\quad(t,u\in F)
\]
with $t^{1+\phi}=u+u^{\phi},$ and the multiplication is given by
\[
Y_{\omega}(t,u)Y_{\omega}(t^{\prime},u^{\prime})=Y_{\omega}(t+t^{\prime
},u+u^{\prime}+t^{\phi}t^{\prime}).
\]

\textbf{Case $B_{2}$:} In this case $p=2,$ $[2]\phi^{2}=1$. The set $\omega$
has type $B_{2}$ with base $\{v,w\}$ where $\langle v,w\rangle=-2$. Elements
of $Y_{\omega}$ take the form
\[
Y_{\omega}(t,u)=X_{v}(t)X_{w}(t^{\phi})X_{v+2w}(u)X_{v+w}(t^{1+\phi}+u^{\phi
})\quad(t,u\in F)
\]
with multiplication $Y_{\omega}(t,u)Y_{\omega}(t^{\prime},u^{\prime
})=Y_{\omega}(t+t^{\prime},u+u^{\prime}+t^{2\phi}t^{\prime})$.\medskip

\textbf{Case $G_{2}$:} Here $p=3,$ $[3]\phi^{2}=1$ and $\left|  F\right|
=3p^{2f}$. The set $\omega$ has type $G_{2}$ with base $\{v,w\}$ where
$\langle v,w\rangle=-3$. Elements of $Y_{\omega}$ take the form
\begin{align*}
Y_{\omega}(t,u,z)  &  =\\
&  X_{v}(t)X_{w}(t^{\phi})X_{v+3w}(u)X_{v+w}(u^{\phi}-t^{1+\phi}%
)X_{2v+3w}(z)X_{v+2w}(z^{\phi}-t^{1+2\phi})
\end{align*}
where $t,\,u,\,z\in F$. The multiplication rule is
\[
Y_{\omega}(t,u,z)Y_{\omega}(t^{\prime},u^{\prime},z^{\prime})=Y_{\omega
}(t+t^{\prime},u+u^{\prime}+t^{\prime}t^{3\phi},z+z^{\prime}-t^{\prime
}u+(t^{\prime})^{2}t^{3\phi}).
\]
\medskip

The root system $\Sigma^{\ast}$ of $S^{\ast}$ is defined as the set of
orthogonal projections $\omega^{\ast}$ of the equivalence classes $\omega$ of
the untwisted root system $\Sigma$ under $\sim$. See definition 2.3.1 of
\cite{GLS} for full details.

The twisted root system $\Sigma^{\ast}$ may not be reduced (i.e. it may
contain several positive scalar multiples of the same root). However in the
case of classical groups $\Sigma^{\ast}$ is reduced with the following
exception: in type $^{2}A_{2m}$ the class $\omega=\{u,v=u^{\tau}, u+v\}$ of
roots in $\Sigma$ spanning a root subsystem of type $A_{2}$ gives rise to a
pair of `doubled' roots $\omega^{\ast}=\{u+v, (u+v)/2\}$ in $\Sigma^{\ast}$.
In this case $\Sigma^{\ast}$ is of type $BC_{m}$, see \cite{GLS} Proposition
2.3.2. Note that the doubled roots $\omega^{\ast}$ above correspond to
\emph{one} root subgroup in $S^{\ast}$, namely $Y_{\omega}$. \medskip

The groups $H,\,U_{+},\,U_{-}\leq S^{\ast}$ are the fixed points of $\sigma$
on the corresponding groups in the untwisted $S$. The group of field
automorphisms $\Phi$ is defined as before; the group of diagonal automorphisms
$\mathcal{D}$ corresponds to the diagonal automorphisms of $S$ that commute
with $\sigma$; there are no graph automorphisms, and we set $\Gamma=1$
\footnote{Note that our definition of graph automorphisms differs from the one
in $\cite{GLS}$}.

\subsection*{The group $\mathcal{D}_{0}$}

In the case when $S$ is a classical group of Lie rank at least 5 we shall
define a certain subgroup $\mathcal{D}_{0}\subseteq\mathcal{D}$ to be used in
Section \ref{smallfield} below.

Suppose first that $S$ is \textbf{untwisted} with a root system $\Sigma$ of
classical type ($A_{r},B_{r},C_{r}$ or $D_{r}$, $r\geq5$) and a set $\Pi$ of
fundamental roots.

If the type is $D_{r}$ define $\Delta=\{w_{1},w_{2}\}$ where $w_{1},w_{2}%
\in\Pi$ are the two roots swapped by the symmetry $\tau$ of $\Pi$. If the type
is $A_{r}$ then let $\Delta:=\{w_{1},w_{2}\}$ where $w_{1},w_{2}\in\Pi$ are
the roots at both ends of the Dynkin diagram (so again we have $w_{1}^{\tau}=
w_{2}$).

If the type is $B_{r}$ or $C_{r}$ set $\Delta=\{w\}$ where $w \in\Pi$ is the
\textbf{long} root at one end of the Dynkin diagram defined by $\Pi$. Recall
that in this case we have $\Gamma=1$.

Let $\Pi_{0}=\Pi\setminus\Delta$ and observe that in all cases $\Pi
_{0}^{\Gamma}=\Pi_{0}$. Now define
\[
\mathcal{D}_{0}= \left\langle h(\chi)\mid\chi_{\left|  \Pi_{0}\right.
}=1\right\rangle .
\]
\medskip

When $S^{\ast}$ is \textbf{twisted} with a root system $\Sigma^{\ast}$ define
$\mathcal{D}_{0}^{\ast} \subseteq S^{\ast}$ to be the group of fixed points of
$\mathcal{D}_{0}$ under $\sigma$, where $\mathcal{D}_{0}$ is the corresponding
subgroup defined for the untwisted version $S$ of $S^{\ast}$. For future
reference, in this case we also consider the set of roots $\Pi_{0}%
\subseteq\Sigma$ as defined above for the untwisted root system $\Sigma$ of
$S$.

When using the notation $S$ for a twisted group, we will write $\mathcal{D}%
_{0}$ for the group here denoted $\mathcal{D}_{0}^{\ast}$.

Clearly $|\mathcal{D}_{0}|\leq|F^{\ast}|^{2},$ and we have

\begin{lemma}
\label{groupD} If $\overline{H}$ is the image of the group $H$ of diagonal
elements in $\mathrm{Inn}(S)$ then $\mathcal{D}=\overline{H}\mathcal{D}_{0}$.

Moreover, provided the type of $S$ is not $A_{r}$ or $^{2}A_{r}$ then there is
a subset $A=\{h(\chi_{i})\ \mid\ 1\leq i\leq4\}\subseteq\mathcal{D}_{0}$ of at
most 4 elements of $\mathcal{D}_{0}$ such that $\mathcal{D}=A\cdot\overline
{H}$.
\end{lemma}

\begin{proof}
We only give the proof for the untwisted case, it easily generalizes to the
twisted case by considering equivalence classes of roots under $\sim$.

From the definition of $\Delta$ one sees that $\Pi_{0}$ can be ordered as
$\nu_{1},\nu_{2},\ldots,\nu_{k}$ so that for some root $w\in\Delta$ we have
\[
\langle\nu_{i},\nu_{i+1}\rangle=\langle\nu_{k},w\rangle=-1,\quad
i=1,2,\ldots,k-1,
\]
and all other possible pairs of roots in $\Pi_{0}\cup\{w\}$ are orthogonal.

Now, given $h=h(\chi)\in\mathcal{D}$ where $\chi$ is a multiplicative
character of the root lattice $\mathbb{Z}\Sigma$, we may recursively define a
sequence $h_{i}=h_{i}(\chi_{i})\in\mathcal{D}$ $(i=1,2,\ldots,k)$ so that
$h_{0}=h$, $h_{i+1}h_{i}^{-1}\in\overline{H}$ and $\chi_{i}$ is trivial on
$\nu_{1},\ldots,\nu_{i}$. Indeed, suppose $h_{i}$ is already defined for some
$i<k$ and $\chi_{i}(\nu_{i+1})=\lambda\in F^{\ast}$, say. Put $h_{i+1}%
=h_{i}h_{\nu_{i+2}}(\lambda)$ (where by convention $\nu_{k+1}=w$). Then
$\chi_{i+1}$ and $\chi_{i}$ agree on $\nu_{1},\ldots,\nu_{i}$, while
\[
\chi_{i+1}(\nu_{i+1})=\chi_{i}(\nu_{i+1})\lambda^{\langle\nu_{i+1},\nu
_{i+2}\rangle}=1.
\]
Clearly we have $h_{k}\in\mathcal{D}_{0}$ while $h\cdot h_{k}^{-1}\in
\overline{H}$. This proves the first statement of the Lemma.

The second statement is now obvious since the group $\mathcal{D}/\overline{H}$
has order at most 4 in that case.
\end{proof}

\medskip

Thus $\mathcal{D}_{0}$ allows us to choose a representative for a given
element in $\mathcal{D}/\overline{H}$ which centralizes many root subgroups
(i.e. those corresponding to $\Pi_{0}$).

\subsection*{Automorphisms and Schur multipliers}

Let $S$ be a Chevalley group as above, untwisted or twisted. We identify
$\overline{S}=S/\mathrm{Z}(S)$ with the group of inner automorphisms
$\mathrm{Inn}(S)$.

\begin{itemize}
\item $\mathrm{Aut}(S)=\overline{S}\mathcal{D}\Phi\Gamma$. (\cite{GLS},
Theorem 2.5.1)

\item $\mathcal{D}\Phi\Gamma$ is a subgroup of $\mathrm{Aut}(S)$ and
$\mathcal{D}\Phi\Gamma\cap\overline{S}=\overline{H}$.

\item When $\Gamma$ is non-trivial it is either of size $2$ or it is
$\mathrm{Sym}(3)$; the latter only occurs in the case $\mathcal{X}=D_{4}$. The
set
\[
\overline{S}\mathcal{D}\Phi
\]
is a normal subgroup of index at most $6$ in $\mathrm{Aut}(S)$.

\item The \emph{universal cover} of $S$ is the largest perfect central
extension $\widetilde{S}$ of $S$. Apart from a finite number of exceptions
$\widetilde{S}$ is the universal Chevalley group of the same type as $S$. The
exceptions arise only over small fields ($\left|  F\right|  \leq9$).

\item The kernel $M(S)$ of the projection $\widetilde{S}\rightarrow S$ is the
\emph{Schur multiplier} of $S$. We have $\left|  M(S)\right|  \leq48$ unless
$S$ is of type $A_{r}$ or $^{2}A_{r}$, in which case (apart from a few small
exceptions) $M(S)$ is cyclic of order dividing $\gcd(r+1,\,\left|
F_{0}\right|  \pm1)$. We also have the crude bound $|M(S)|\leq|S|$.

\item $\left|  \mathcal{D}:\overline{H}\right|  \leq\left|  M(\overline
{S})\right|  $

\item $\left|  \mathrm{Out}(S)\right|  \leq2f|M(\overline{S})|$ where $\left|
F\right|  =p^{f}$ unless $S$ is of type $D_{4}$.
\end{itemize}

Suppose that $T$ is a quasisimple group of Lie type, with $T/\mathrm{Z}(T)=S$.
Then $T=\widetilde{S}/K$ for some $K\leq\mathrm{Z}(\widetilde{S})$. An
automorphism $\gamma$ of $\widetilde{S}$ that stabilizes $K$ induces an
automorphism $\overline{\gamma}$ of $T$. The map $\gamma\mapsto\overline
{\gamma}$ is an isomorphism between $N_{\mathrm{Aut}(\widetilde{S})}(K)$ and
$\mathrm{Aut}(T)$; see \cite{A}, section 33. Thus every automorphism of $T$
lifts to an automorphism of $\widetilde{S}.$

\section{Combinatorial lemmas}

The first three lemmas are elementary, and we record them here for
convenience. $G$ denotes an arbitrary finite group.

\begin{lemma}
\label{m-prod}Suppose that $\left|  G\right|  \leq m$.

\emph{(i)} If $f_{1},\ldots,f_{m}\in G$ then $\prod_{l=i}^{j}f_{l}=1$ for some
$i\leq\,j$.

\emph{(ii)} If $G=\left\langle X\right\rangle $ and $1\in X$ then $G=X^{\ast
m}.$
\end{lemma}

\begin{lemma}
\label{trick} Let $M$ be a $G$-module and suppose that $\sum_{i=1}^{L}%
M(g_{i}^{e_{i}}-1)=M$ for some $g_{i}\in G$ and $e_{i}\in\mathbb{N}$. Then
\[
M=\sum_{i=1}^{L}M(g_{i}-1).
\]

\end{lemma}

\begin{lemma}
\label{pr} Let $\alpha_{1},\alpha_{2},\ldots,\alpha_{m}\in\mathrm{Aut}(G)$.
Then
\[
\lbrack G,\alpha_{1}\alpha_{2}\ldots\alpha_{m}]\subseteq\lbrack G,\alpha
_{1}]\cdot\ldots\cdot\lbrack G,\alpha_{m}].
\]

\end{lemma}

We shall also need the following useful result, due to Hamidoune:

\begin{lemma}
\label{ham}\emph{\cite{H}} Let $X$ be a subset of $G$ such that $X$ generates
$G$ and $1\in X$. If $|G|\leq m\left|  X\right|  $ then $G=X^{\ast2m}$.
\end{lemma}

We conclude with some remarks about quasisimple groups. Let $S$ be a finite
quasisimple group. Then $\mathrm{Aut}(S)$ maps injectively into $\mathrm{Aut}%
(\overline{S})$ and $\mathrm{Out}(S)$ maps injectively into $\mathrm{Out}%
(\overline{S})$. Since every finite simple group can be generated by $2$
elements \cite{AG}, it follows that
\[
|\mathrm{Aut}(S)|\leq|\mathrm{Aut}(\overline{S})|\leq\left|  \overline
{S}\right|  ^{2}.
\]
Also $S$ can be generated by $2$ elements, since if $S=\left\langle
X\right\rangle \mathrm{Z}(S)$ then $S=[S,S]\subseteq\left\langle
X\right\rangle $.

Since $\left|  M(\overline{S})\right|  <\left|  \overline{S}\right|  $
(\cite{G}, table 4.1\textbf{)} and $\left|  \mathrm{Z}(S)\right|  \leq\left|
M(\overline{S})\right|  $ we have $\left|  S\right|  \leq\left|  \overline
{S}\right|  ^{2}$.

If $g\in S\setminus\mathrm{Z}(S)$ then $[S,g]\cdot\lbrack S,g]^{-1}$ contains
(many) non-central conjugacy classes of $S$; it follows that $S$ is generated
by the set $[S,g]$.

\section{Deduction of Theorem \ref{twisted}\label{deduction}}

This depends on the special case of Theorem \ref{autos} where $q=1$. Assuming
that this case has been proved, we begin by showing that the constant $C(1)$
may be reduced to $1$, provided the constant $M(1)$ is suitably enlarged.

Let $\mathcal{S}$ denote the finite set of quasisimple groups $S$ such that
$|S/\mathrm{Z}(S)|\leq C=C(1)$, and put $M^{\prime}=C^{4}.$ We claim that if
$S\in\mathcal{S}$ and $\beta_{1},\ldots,\beta_{M^{\prime}}$ are any
automorphisms of $S$ then there exist $g_{1},\ldots,g_{M^{\prime}}\in S$ such
that
\begin{equation}
S=\prod_{j=1}^{M^{\prime}}[S,\overline{g}_{j}\beta_{j}]. \label{prod-comm}%
\end{equation}
Thus in Theorem \ref{autos} we may replace $C(1)$ by $1$ provided we replace
$M(1)$ by $\max\{M(1),\,M^{\prime}\}$.

Since $|\mathrm{Aut}(S)|\leq\left|  \overline{S}\right|  ^{2}\leq C^{2}$,
Lemma \ref{m-prod}(i) implies that the sequence \linebreak$(\beta_{1}%
,\ldots,\beta_{M^{\prime}})$ contains subsequences $(\beta_{1}(i),\ldots
,\beta_{j(i)}(i))$, $i=1,\ldots,C^{2}$, such that (a) $\prod_{l=1}^{j(i)}$
$\beta_{l}(i)=1$ for each $i$, and (b) for each $i<C^{2} $, $\beta_{j(i)}(i)$
precedes $\beta_{1}(i+1)$; we will call such subsequences `strictly disjoint'.

Fix a non-central element $g\in S$ and put
\begin{align*}
g_{1}(i)  &  =g\\
g_{2}(i)  &  =\cdots=g_{j(i)}(i)=1
\end{align*}
for $i=1,\ldots,C^{2}$. Then Lemma \ref{pr} gives
\[
\prod_{l=1}^{j(i)}[S,\overline{g_{l}(i)}\beta_{l}(i)]\supseteq\lbrack
S,\overline{g}\prod_{l=1}^{j(i)}\beta_{l}(i)]=[S,g]
\]
for each $i$. As $[S,g]$ generates $S$ and $\left|  S\right|  <\left|
\overline{S}\right|  ^{2}\leq C^{2}$ it now follows by Lemma \ref{m-prod}(ii)
that
\[
S=[S,g]^{\ast C^{2}}\subseteq\prod_{i=1}^{C^{2}}\prod_{l=1}^{j(i)}%
[S,\overline{g_{l}(i)}\beta_{l}(i)].
\]
This gives (\ref{prod-comm}) for a suitable choice of the $g_{j}$, each equal
to either $g$ or $1$.

\bigskip

Let us re-define $M(1)$ now so that the statement of Theorem \ref{autos} holds
with $C(1)=1$, and set $D=M(1).$ Then Theorem \ref{twisted} follows on taking
$S=G$ and $k=D$ in the next lemma:

\begin{lemma}
\label{untwisted}Let $G$ be a group and let $\beta_{1},\ldots,\,\beta_{k}%
\in\mathrm{Aut}(G)$. Suppose that there exist $g_{1},\ldots,g_{k}\in G$ such
that
\[
G=\prod_{i=1}^{k}[G,\overline{g}_{i}\beta_{i}].
\]
Then for any $\alpha_{1},\ldots,\alpha_{k}\in\mathrm{Aut}(G)$ we have
\[
G=\prod_{i=1}^{k}T_{\alpha_{i},\beta_{i}}(G,G).
\]

\end{lemma}

\begin{proof}
Note the identities
\begin{align*}
T_{\alpha,\beta}(h^{-\alpha^{-1}},z^{h})  &  =[h^{-1},\alpha^{-1}]^{-1}%
\cdot\lbrack z,\overline{h}\beta],\\
a[xa,\beta]  &  =[x,\beta]a^{\beta}\text{;}%
\end{align*}
the second one implies that
\[
a[G,\beta]=[G,\beta]a^{\beta}%
\]
for any $a\in G$ and $\beta\in\mathrm{Aut}(G)$.

Now let $g_{1},\ldots,g_{k}$ be the given elements of $G$ and put
$a_{i}=[g_{i}^{-1},\alpha_{i}^{-1}]^{-1}$ for each $i$. Then
\begin{align*}
\prod_{i=1}^{k}T_{\alpha_{i},\beta_{i}}(G,G)\supseteq a_{1}[G,\overline{g}%
_{1}\beta_{1}]a_{2}[G,\overline{g}_{2}\beta_{2}]\ldots a_{k}[G,\overline
{g}_{k}\beta_{k}]\\
=[G,\overline{g}_{1}\beta_{1}][G,\overline{g}_{2}\beta_{2}]\ldots\lbrack
G,\overline{g}_{k}\beta_{k}]\cdot b=G
\end{align*}
where
\[
b=\prod_{i=1}^{k}a_{i}^{\overline{g}_{i}\beta_{i}\ldots\overline{g}_{k}%
\beta_{k}}.
\]
This completes the proof.
\end{proof}

\section{Alternating groups and groups of Lie type over small fields}

\label{smallfield}

Given $q\in\mathbb{N}$ we fix a large integer $K=K(q)$ (greater than $100$,
say); how large $K$ has to be will appear in due course. Let $\mathcal{S}%
_{1\mathbf{a}}$ denote the family of all quasisimple groups $S$ such that
$\overline{S}=\mathrm{Alt}(k)$ for some $k>K$, and let $\mathcal{S}%
_{1\mathbf{b}}$ denote the family of all quasisimple groups of Lie type of Lie
rank greater than $K$ over fields $F$ with $\left|  F\right|  \leq K$. Let
$\mathcal{S}_{1}=\mathcal{S}_{1\mathbf{a}}\cup\mathcal{S}_{1\mathbf{b}}$.

For $S\in\mathcal{S}_{1}$ we define a subgroup $S_{0}$ of $S$ as follows:

\begin{description}
\item[ (a)] if $\overline{S}=\mathrm{Alt}(k)$, let $S_{0}$ be the inverse
image in $S$ of $\mathrm{Alt}(\left\{  3,\ldots,k\right\}  )\cong
\mathrm{Alt}(k-2)$, the pointwise stabilizer of $\{1,2\}$;

\item[ (b)] if $S$ is of Lie type, first recall the definition of
$\mathcal{D}_{0}\subseteq\mathcal{D}$ and $\Pi_{0}\subseteq\Sigma$ from
Section \ref{resume}.

In case $\Sigma=A_{r},D_{r}$ (i.e. $S$ has type $A_{r},D_{r}$, $^{2}A_{r}$ or
$^{2}D_{r}$) there is a graph automorphism $\tau\not =1$ of the untwisted
version of $S$. Define $S_{0}$ to be the group of fixed points under $\tau$ of
the group
\[
R:=\left\langle X_{v}(\lambda)\mid v\in\pm\Pi_{0}, \ \lambda\in\mathbb{F}%
_{p}\right\rangle .
\]
Here $\mathbb{F}_{p}$ is the prime field of $F$.

In the remaining cases (i.e. $\Sigma$ of type $B_{r},C_{r}$ and $S$ is
untwisted) define $S_{0}:=R$.

It is easy to see that in all cases from (\textbf{b}) we have $S_{0}\leq S$
and $S_{0}$ is centralized by $\mathcal{D}_{0}\Phi\Gamma$.
\end{description}

\bigskip

In case (\textbf{a}), let $\tau$ denote the lift to $\mathrm{Aut}(S)$ of the
automorphism of $\overline{S}$ given by conjugation by $(12)$. Then
$\mathrm{Aut}(S)=\mathrm{Inn}(S)\left\langle \tau\right\rangle $ and $\tau$
acts trivially on $S_{0}$. Also $\log|S|/\log|S_{0}|\leq2$ and $\left|
\mathrm{Z}(S)\right|  \leq2$.

In case (\textbf{b}), $S_{0}$ is again a quasisimple group of Lie type, and of
Lie rank at least $K/2-1$ (\cite{GLS}, \S 2.3). It is fixed elementwise by
automorphisms of $S$ lying in the set $\mathcal{D}_{0} \Phi\Gamma$, and we
have
\begin{equation}
\log|S|/\log|S_{0}|\leq2(F:\mathbb{F}_{p})+A\leq3\log K, \label{SandS0}%
\end{equation}
say, where $A$ is some constant. Also $\left|  \mathrm{Z}(S)\right|  \leq K$.

From Section \ref{resume} we have $|\mathcal{D}_{0}|\leq\left|  F^{\ast
}\right|  ^{2}$ and $\mathrm{Aut}(S)=\mathrm{Inn}(S)\mathcal{D}_{0}\Phi\Gamma
$\textbf{.} Note that $\mathcal{D}_{0}\Phi\Gamma$ is a subgroup of
$\mathrm{Aut}(S)$ and $\left|  \mathcal{D}_{0}\Phi\Gamma\right|  \leq
2K^{2}\log K$.

Now let $S\in\mathcal{S}_{1}$. To prove Theorem \ref{autos} for $S$, we have
to show that given automorphisms $\beta_{1},\ldots,\beta_{M}$ of $S$, where
$M$ is sufficiently large, and given divisors $q_{1},\ldots,q_{M}$ of $q$, we
can find inner automorphisms $\alpha_{1},\ldots,\alpha_{M}$ of $S$ such that
\begin{equation}
S=[S,(\alpha_{1}\beta_{1})^{q_{1}}]\cdot\ldots\cdot\lbrack S,(\alpha_{M}%
\beta_{M})^{q_{M}}]. \label{statement}%
\end{equation}
We may freely adjust each of the automorphisms $\beta_{i}$ by an inner
automorphism, so without loss of generality we assume that each $\beta_{i}%
\in\{1,\tau\}$ if $\overline{S}$ is alternating, and that each $\beta_{i}%
\in\mathcal{D}_{0}\Phi\Gamma$ if $S$ is of Lie type.

\begin{lemma}
\label{q} Provided $K=K(q)$ is sufficiently large, there exists a $q$th power
$h$ in $S_{0}$ such that
\[
\log|\overline{h}^{\overline{S}}|\geq\frac{\log|S|}{36\log K}.
\]

\end{lemma}

\begin{proof}
By examining the proofs of Lemmas 1 -- 4 of \cite{SW}, one finds that provided
the degree (resp the Lie rank) of $\overline{S_{0}}$ is large enough,
$\overline{S_{0}}$ is a product of six conjugacy classes of $q$th powers (even
stronger results may be deduced from \cite{LS2}.) It follows that $S_{0}$
contains a $q$th power $h$ such that the conjugacy class of $\overline{h}$ in
$\overline{S_{0}}$ has size at least $|\overline{S_{0}}|^{1/6}$. The result
now follows from (\ref{SandS0}) since $|S_{0}|\leq|\overline{S_{0}}|^{2}$.
\end{proof}

\medskip

Next, we quote a related result of Liebeck and Shalev:

\begin{theorem}
\emph{(\cite{LS2}, Theorem 1.1) }There is an absolute constant $C_{0}$ such
that for every simple group $T$ and conjugacy class $X$ of $T$ we have
\[
T=X^{\ast t},
\]
where $t=\left\lfloor C_{0}\log|T|/\log|X|\right\rfloor $.
\end{theorem}

Now take $T=S/\mathrm{Z}(S)$ and $X=\overline{h}^{T}$ where $h\in S_{0}$ is
given by Lemma \ref{q}. Since $[S,h]=(h^{-1})^{S}\cdot h$, the above theorem
gives
\[
S=[S,h]^{\ast\left\lceil 36C_{0}\log K\right\rceil }\cdot\mathrm{Z}(S),
\]
and applying Lemma \ref{ham} we deduce that
\begin{equation}
S=[S,h]^{\ast M^{\prime}} \label{[S,h]}%
\end{equation}
where $M^{\prime}=\left\lceil 72KC_{0}\log K\right\rceil $.

Suppose now that $M\geq2K^{2}\log K\cdot M^{\prime}$. Then the group generated
by $\beta_{1},\ldots,\beta_{M}$ has order at most $M/M^{\prime},$ and we may
use Lemma \ref{m-prod}, as in the preceding section, to find strictly disjoint
subsequences $(\beta_{1}(i)^{q_{1,i}},\ldots,\beta_{j(i)}(i)^{q_{j(i),i}})$ of
$(\beta_{1}^{q_{1}},\ldots,\beta_{M}^{q_{M}}) $, $i=1,\ldots,M^{\prime}$, such
that $\prod_{l=1}^{j(i)}\beta_{l}(i)^{q_{l,i}}=1$ for each $i$.

Since $h$ is is a $q$th power in $S_{0}$, for each $i$ there exists $h_{i}\in
S_{0}$ such that $h_{i}^{q_{1,i}}=h$. Then each $\overline{h_{i}}$ commutes
with each $\beta_{j}$; so putting
\begin{align*}
\alpha_{1}(i)  &  =\overline{h_{i}}\\
\alpha_{j}(i)  &  =1\qquad(j\geq2)
\end{align*}
we have
\[
\prod_{l=1}^{j(i)}(\alpha_{l}(i)\beta_{l}(i))^{q_{l,i}}=\overline{h}.
\]
Hence
\[
\lbrack S,h]=[S,\prod_{l=1}^{j(i)}(\alpha_{l}(i)\beta_{l}(i))^{q_{l,i}%
}]\subseteq\prod_{l=1}^{j(i)}[S,(\alpha_{l}(i)\beta_{l}(i))^{q_{l,i}}]
\]
by Lemma \ref{pr}, and (\ref{statement}) follows from (\ref{[S,h]}).

Thus Theorem \ref{autos} holds for groups $S\in\mathcal{S}_{1}$ provided
$K=K(q)$ and $M=M(q)$ are sufficiently large.

\section{Groups of Lie type over large fields: reductions}

\label{Orb}

As before, we fix a positive integer $q$ and denote by $K=K(q)$ some large
positive integer, to be specified later. Let $\mathcal{S}_{2}$ (resp.
$\mathcal{S}_{3}$) denote the family of all quasisimple groups of Lie type of
Lie rank at most $8$ (resp. at least $9$) over fields $F$ with $\left|
F\right|  >K$. According to the Classification, all but finitely many finite
quasisimple groups lie in $\mathcal{S}_{1}\cup\mathcal{S}_{2}\cup
\mathcal{S}_{3}$, so it remains only to prove Theorem \ref{autos} for groups
in $\mathcal{S}_{2}\cup\mathcal{S}_{3}$.

The validity of this theorem for groups of Lie type over large fields depends
on there being `enough room' for certain equations to be solvable. In order to
exploit this, we need to restrict the action of the relevant automorphisms to
some very small subgroups; this is made possible (as in the preceding section)
by choosing a suitable representative for each outer automorphism. The desired
`global' conclusion will then be derived with the help of the following result
of Liebeck and Pyber:

\begin{theorem}
\emph{(\cite{LP}, Theorem D) }Let $S$ be a quasisimple group of Lie type. Then
$S=(U_{+}U_{-})^{\ast12}\cdot U_{+}$.
\end{theorem}

For the rest of this section, $S$ denotes a group in $\mathcal{S}_{2}%
\cup\mathcal{S}_{3}$, and we use the notation of Section \ref{resume} for root
subgroups etc. Our aim is to prove

\begin{proposition}
\label{uni}There is a constant $M_{1}=M_{1}(q)$ such that if $\gamma
_{1},\ldots,\gamma_{M_{1}}$ are automorphisms of $S$ lying in $\mathcal{D}%
\Phi\Gamma$ and $q_{1},\ldots,q_{M_{1}}$ are divisors of $q$ then there exist
elements $h_{1},\ldots,h_{M_{1}}\in H$ and $u_{1},\ldots,u_{M_{1}}\in U$ such
that
\[
U\subseteq\prod_{i=1}^{M_{1}}[U,(\overline{u_{i}h_{i}}\gamma_{i})^{q_{i}}].
\]

\end{proposition}

By symmetry, the same result will then hold with $U_{-}$ in place of $U=U_{+}
$. Since $\mathrm{Aut}(S)=\mathrm{Inn}(S)\mathcal{D}\Phi\Gamma$, this together
with the above theorem shows that Theorem \ref{autos} holds for $S$ as long as
$M(q)\geq25M_{1}(q).$

To establish Proposition \ref{uni}, we shall express $U$ as a product of
certain special subgroups, each of which itself satisfies a similar property.
There are several different cases to consider, and we apologise for the
complexity of the argument. The basic idea in all cases is the same: the
required result is reduced to showing that certain equations are solvable over
a suitable finite field, and then applying a general result about such
equations, namely Lemma \ref{sur} below.

The first class of special subgroups is defined as follows:

\begin{definition}
Let $S$ be a quasisimple group of Lie type.\newline\medskip\textbf{Case 1:}
when $S$ is \emph{untwisted}\textbf{. }Let $\omega$ be an equivalence class of
roots from $\Sigma/\sim$ as defined in Section \ref{resume}. The corresponding
\emph{orbital subgroup} is then defined to be
\[
O(\omega)=W_{\omega}:=\left\langle X_{v}\mid v\in\omega\right\rangle
=\prod_{v\in\omega}X_{v}%
\]
(product ordered by increasing height of roots).\newline\medskip\textbf{Case
2:} when $S=S^{\ast}$ is \emph{twisted}. Define
\[
O(\omega):=Y_{\omega}%
\]
to be the root subgroup of $S$ corresponding to the equivalence class $\omega$
of the untwisted root system $\Sigma$ under $\sim$ described in Section
\ref{resume}.
\end{definition}

\noindent Note (i) the orbital subgroups are invariant under $\mathcal{D}%
\Phi\Gamma$. (ii) In case 2, $Y_{\omega}$ is in fact a subgroup of the orbital
subgroup $W_{\omega}$ defined in Case 1 for the untwisted version of $S^{\ast
}$.

In Section \ref{orbitalsec} we prove

\begin{proposition}
\label{orb}There is a positive integer $L=L(q)$ such that if $\gamma
_{1},\ldots,\gamma_{L}$ are automorphisms of $S$ lying in $\mathcal{D}%
\Phi\Gamma$, $q_{1},\ldots,q_{L}$ are divisors of $q$ and $O=O(\omega)$ is an
orbital subgroup of $S$ then there exist elements $h_{1},\ldots,h_{L}\in H$
such that
\[
O\subseteq\prod_{i=1}^{L}[O,(\overline{h_{i}}\gamma_{i})^{q_{i}}].
\]

\end{proposition}

Now suppose that $S\in\mathcal{S}_{2}$. Then $U$ is equal to the product of at
most $120$ root subgroups, each of which is contained in an orbital subgroup.
So in this case Proposition \ref{uni} follows as long as we take $M_{1}%
(q)\geq120L$.

\bigskip

For arbitrary groups $S$ $\in\mathcal{S}_{3}$ we can't write $U$ as a bounded
product of orbital subgroups. However, $S$ is then a classical group, and
contains a relatively large subgroup of type $\mathrm{SL}$. For the group
$\mathrm{SL}$ itself we prove

\begin{proposition}
\label{SL}Let $S=\mathrm{SL}_{r+1}(F)$, where $\left|  F\right|  >K$ and
$r\geq3$. There is a constant $M_{2}=M_{2}(q)$ such that if $\gamma_{1}%
,\ldots,\gamma_{M_{2}}$ are automorphisms of $S$ lying in $\mathcal{D}%
\Phi\Gamma$ and $q_{1},\ldots,q_{M_{2}}$ are divisors of $q$ then there exist
automorphisms $\eta_{1},\ldots,\eta_{M_{2}}\in\mathcal{D}$ and elements
$u_{1},\ldots,u_{M_{2}}\in U$ such that
\[
U\subseteq\prod_{i=1}^{M_{2}}[U,(\overline{u_{i}}\eta_{i}\gamma_{i})^{q_{i}%
}].
\]

\end{proposition}

\noindent Note that this differs from Proposition \ref{uni} in that the
elements $\eta_{i}$ are allowed to vary over \emph{all diagonal automorphisms
}of $S$, not just the inner ones.

We also need to consider the following special subgroups $V_{s+1}$ and
$V_{s+1}^{\ast}$ of the full unitriangular group:

\begin{definition}
\emph{(i)} In $\mathrm{SL}_{s+1}(F)$ define%
\[
V_{s+1}=\left(
\begin{array}
[c]{ccccc}%
1 & \ast & \cdots & \ast & \ast\\
& 1 & \mathbf{0} & 0 & \ast\\
&  & \ddots & \mathbf{0} & \vdots\\
&  &  & 1 & \ast\\
&  &  &  & 1
\end{array}
\right)  <\mathrm{SL}_{s+1}(F).
\]
the group of unitriangular matrices differing from the identity only in the
first row and last column.\newline\emph{(ii) }Consider $\mathrm{SU}_{s+1}(F)$
as the set of fixed points of the Steinberg automorphism $\sigma$ of
$\mathrm{SL}_{s+1}(F)$. Then%
\[
V_{s+1}^{\ast}=V_{s+1}\cap\mathrm{SU}_{s+1}(F),
\]
the set of fixed points of $\sigma$ on $V_{s+1}$.
\end{definition}

\noindent Note that $V_{s+1}$ and $V_{s+1}^{\ast}$ are stabilized by
automorphisms of $\mathrm{SL}_{s+1}(F)$, resp. $\mathrm{SU}_{s+1}(F)$, lying
in $\mathcal{D}\Phi\Gamma$.

\begin{proposition}
\label{manyW}Let $S=\mathrm{SL}_{s+1}(F)$ or $\mathrm{SU}_{s+1}(F)$ where
$s\geq5$ and $\left|  F\right|  >K$. Put $V=V_{s+1}$ in the first case,
$V=V_{s+1}^{\ast}$ in the second case. There is a positive integer
$L_{1}=L_{1}(q)$ such that if $\gamma_{1},\ldots,\gamma_{L_{1}}$ are
automorphisms of $S$ lying in $\mathcal{D}\Phi\Gamma$ and $q_{1}%
,\ldots,q_{L_{1}}$ are divisors of $q$ then there exist elements $h_{1}%
,\ldots,h_{L_{1}}\in H$ such that
\[
V=\prod_{i=1}^{L_{1}}[V,(\overline{h_{i}}\gamma_{i})^{q_{i}}].
\]

\end{proposition}

\noindent Of course, the point here is that $L_{1}$ is independent of $s$.

The last two propositions will be proved in Section \ref{unitri}. We need one
more kind of special subgroup, denoted $P$. This is defined for groups $S$ of
type%
\[
\mathcal{X}\in\{^{2}A_{r},\,B_{r},\,C_{r},\,D_{r},\,^{2}D_{r}\},\,\,r\geq4.
\]
Recall that $\Sigma$ is the root system of $\mathcal{X}$ (twisted or
untwisted). In each of these five cases, there exist fundamental roots
$\delta,\,\delta^{\prime}\in\Sigma$ (equal unless $\mathcal{X}=D_{r}$, see
below) such that the other fundamental roots $\Pi^{\prime}=\Pi-\{\delta
,\delta^{\prime}\}$ generate a root system $\Sigma^{\prime}$ of type $A_{s}$,
for the appropriate $s\,$($=\left\lceil \frac{r}{2}\right\rceil -1$, $\ r-1$
or $r-2$): in types $^{2}A_{r},B_{r},C_{r}$ and $^{2}D_{r}$ we take
$\delta=\delta^{\prime}$ to be the fundamental root of length distinct from
the others; in type $D_{r}$, $\{\delta,\delta^{\prime}\}$ is the pair of
fundamental roots swapped by the symmetry $\tau$ of $D_{r}$ (see \cite{GLS},
Proposition 2.3.2).

If $S$ is untwisted, put $S_{1}=\left\langle X_{w}\mid w\in\Sigma^{\prime
}\right\rangle $ and $U_{1}=\prod_{w\in\Sigma_{+}^{\prime}}X_{w}$ . If $S$ is
twisted, define $S_{1}$ and $U_{1}$ similarly by replacing the root subgroups
$X_{w}$ with the corresponding root subgroup $Y_{\omega}$, $\omega^{\ast}%
\in\Sigma^{\prime}$.

Then $S_{1}$ is a quasisimple group of type $A_{s}$ (a Levi subgroup of $S$),
it is fixed pointwise by $\Gamma$, and $U_{1}$ is its positive unipotent subgroup.

\begin{definition}
\label{Pdef}If $S$ is untwisted, set%
\[
P=\prod_{w\in\Sigma_{+}\backslash\Sigma_{+}^{\prime}}X_{w}.
\]
If $S$ is twisted, set%
\[
P=\prod_{\omega^{\ast}\in\Sigma_{+}\backslash\Sigma_{+}^{\prime}}Y_{\omega}.
\]

\end{definition}

\noindent Note that $P$ is a subgroup of $U$ stabilized by $\mathcal{D}%
\Phi\Gamma$ and that
\[
U=U_{1}P.
\]
In the final section we prove

\begin{proposition}
\label{P1}Assume that $\left|  F\right|  >K$ and that $S$ is of type
$B_{r},\,C_{r},\,D_{r}$ or $\,^{2}D_{r}$, where $r\geq4$. There is a constant
$N_{1}=N_{1}(q)$ such that if $\gamma_{1},\ldots,\gamma_{N_{1}}$ are
automorphisms of $S$ lying in $\mathcal{D}\Phi\Gamma$ and $q_{1}%
,\ldots,q_{N_{1}}$ are divisors of $q$ then there exist elements $h_{1}%
,\ldots,h_{N_{1}}\in H$ such that%
\[
P\subseteq\prod_{i=1}^{N_{1}}[P,(\overline{h_{i}}\gamma_{i})^{q_{i}}].
\]

\end{proposition}

\begin{proposition}
\label{P2}Assume that $\left|  F\right|  >K$ and that $S$ is of type
$\,^{2}A_{r}$, where $r\geq4$. There is a constant $N_{1}^{\prime}%
=N_{1}^{\prime}(q)$ such that if $\gamma_{1},\ldots,\gamma_{N_{1}^{\prime}}$
are automorphisms of $S$ lying in $\mathcal{D}\Phi\Gamma$ and $q_{1}%
,\ldots,q_{N_{1}^{\prime}}$ are divisors of $q$ then there exist automorphisms
$\eta_{1},\ldots,\eta_{N_{1}^{\prime}}\in\mathcal{D}$ such that
\[
P\subseteq\prod_{i=1}^{N_{1}^{\prime}}[P,(\eta_{i}\gamma_{i})^{q_{i}}].
\]

\end{proposition}

Assuming the last five propositions, we may now complete the\medskip

\noindent\textbf{Proof of Proposition \ref{uni}.} We take\textbf{ }$S$ to be a
quasisimple group in $\mathcal{S}_{3}$, and deal with each type in
turn.\medskip

\textbf{Case 1.} \ Type $A_{s}$, where $s\geq5$. (For $S\in\mathcal{S}_{3}$ we
only need $s\geq9$; the more general result is needed for later applications.)

Assume to begin with that $S=\mathrm{SL}_{s+1}(F)$; here $\left|  F\right| >K$
and $s\geq5$. Let $S^{1}$ be the copy of $\mathrm{SL}_{s-1}(F)$ sitting in the
middle $(s-1)\times(s-1)$ square of $S$, and let $U^{1}$ denote the upper
unitriangular subgroup of $S^{1}$.

Now $U^{1}$ is the positive unipotent subgroup of $S^{1}$, and the diagonal
subgroup $H$ of $S$ induces by conjugation on $S^{1}$ its full group of
diagonal automorphisms. Also, diagonal, field and graph automorphisms of $S$
restrict to automorphisms of the same type on $S^{1}$. Thus given $M_{2}$
automorphisms $\gamma_{i}$ of $S$ lying in $\mathcal{D}\Phi\Gamma$ and $M_{2}$
divisors $q_{i}$ of $q$, Proposition \ref{SL} applied to $S^{1}$ shows that
there exist elements $h_{i}\in H$ and $u_{i}\in U^{1}$ such that%

\[
U^{1}\subseteq\prod_{i=1}^{M_{2}}[U^{1},(\overline{u_{i}h_{i}}\gamma
_{i})^{q_{i}}].
\]

On the other hand, we have $U=U^{1}V_{s+1}$. With Proposition \ref{manyW} this
shows that Proposition \ref{uni} holds for $S=\mathrm{SL}_{s+1}(F)$ provided
we take $M_{1}\geq M_{2}+L_{1}$.

It then follows for every quasisimple group $S$ of type $A_{s}$ ($s\geq5$),
since automorphisms of $S$ in $\mathcal{D}\Phi\Gamma$ lift to automorphisms of
the same type of its covering group $\mathrm{SL}_{s+1}(F)$.

\medskip

\textbf{Case 2.} Type $^{2}\!A_{s}$, where $s\geq9$.

As above, we may assume that in fact $S=\mathrm{SU}_{s+1}(F)$. Considering the
fixed points of the Steinberg automorphism in $\mathrm{SL}_{s+1}(F)$, we see
that
\[
U=U^{1}V_{s+1}^{\ast}%
\]
where $U^{1}$ is the positive unipotent subgroup of $S^{1}$ and $S^{1}$ is a
copy of $\mathrm{SU}_{s-1}(F)$ sitting `in the middle' of $S$. Again, the
diagonal subgroup $H$ of $S$ induces by conjugation on $S^{1}$ its full group
of diagonal automorphisms.

Now we apply Definition \ref{Pdef} and the discussion preceding it to the
group $S^{1}$; this gives a subgroup $S_{2}$ of $S^{1}$ of type $A_{t}$ where
$t=\left\lceil \frac{s-2}{2}\right\rceil -1\geq3$, and the subgroup $P$, and
we have
\[
U^{1}=U_{2}P
\]
where $U_{2}$ denotes the positive unipotent subgroup of $S_{2}$.

Proposition \ref{uni} now follows for $S$ on combining Propositions \ref{SL}
(for $U_{2}$), \ref{P2} (for $P$), and \ref{manyW} (for $V_{s+1}^{\ast}$), and
taking $M_{1}=M_{2}+L_{1}+N_{1}^{\prime}$.

\medskip

\textbf{Case 3.} Type $B_{r},\,C_{r},\,D_{r}$ or $\,^{2}D_{r}$,
where$\,\,r\geq9$.

Again we apply Definition \ref{Pdef} and the discussion preceding it. This
gives a subgroup $S_{1}$ of $S$ of type $A_{s}$, where $s\geq7$, and the
corresponding subgroup $P$, and $U=U_{1}P$ where $U_{1}$ is the positive
unipotent subgroup of $S_{1}$. In this case, Proposition \ref{uni} follows
from Case 1, done above, and Proposition \ref{P1}, on taking $M_{1}%
=M_{2}+L_{1}+N_{1}$.

\medskip This completes the proof, modulo Propositions 6.4 -- 6.10.

The proofs of Propositions 6.9 and 6.10, given in Section 10, are similar to
that of Proposition 6.4, but even more complicated. An alternative approach is
available if we are prepared to quote the following general result, which will
be established in \cite{N}:

\begin{theorem}
Let $S$ be a quasisimple group of classical type over a finite field $F$. Then
$S$ contains a subgroup $S_{1}$ such that\newline\emph{(i) }$S_{1}$ is
quasisimple of type $A_{n}$ over $F$ or $F_{0}$, for some $n$,\newline%
\emph{(ii)} $S_{1}$ is invariant under $\mathcal{D}\Phi\Gamma,$\newline%
\emph{(iii)} there exist $g_{1},\ldots,g_{h}\in S$ such that
$S={\displaystyle\prod\limits_{i=1}^{h}}g_{i}^{-1}S_{1}g_{i}$, where $h$ is an
absolute constant.
\end{theorem}

\noindent(The constant $h$ is probably rather less than $600$.)

Suppose we have established Proposition 6.2 for groups of type $A_{n}$. Then
Theorem 1.2 holds for such groups, with a constant $M(q)$ say. Now let
$S\in\mathcal{S}_{3}$, not of type $A_{n}$, and let $\beta_{ij}$ ($1\leq i\leq
h,~~1\leq j\leq M(q)$) be given automorphisms of $S$. Let $x_{ij}\in S$ be
such that $\gamma_{ij}:=\overline{x_{ij}}\beta_{ij}\in\mathcal{D}\Phi\Gamma$.
Applying Theorem 1.2 to $S_{1}$ we find elements $y_{ij}\in S_{1}$ such that
for each $i$,
\[
S_{1}={\displaystyle\prod\limits_{j=1}^{M(q)}}[S_{1},~(\overline{y_{ij}}%
\gamma_{ij})^{q_{ij}}].
\]
With (iii) this gives
\[
S={\displaystyle\prod\limits_{i=1}^{h}}{\displaystyle\prod\limits_{j=1}%
^{M(q)}} [S_{1}^{g_{i}},~(\overline{u_{ij}}\beta_{ij})^{q_{ij}}]
\]
where
\[
u_{ij}=g_{i}^{-1}y_{ij}x_{ij}g_{i}^{\beta_{ij}^{-1}}\in S.
\]

Thus Theorem 1.2 holds for $S$ on replacing $M(q)$ by $hM(q)$. In this
approach, it suffices to prove only Case 1 of Proposition 6.2, and Definition
6.8 and Propositions 6.9, 6.10 may be omitted (as may the part of Proposition
6.7 dealing with the group $\mathrm{SU}$).

\section{Surjective maps}

The following lemma lies at the heart of the proof; here $q$ and $M$ denote
positive integers:

\begin{lemma}
\label{sur} Let $F$ be a finite field and let $\phi_{1},\ldots,\phi_{M}$ be
automorphisms of $F$. Let $\mu_{1},\mu_{2},\ldots,\mu_{M}$ be nonzero elements
of $F$, let $q_{1},\ldots,q_{M}$ be divisors of $q$ and let $c_{1}%
,\ldots,c_{M}$ be positive integers. For $\lambda,\,t\in F$ put
\[
f_{i,\lambda}(t)=\mu_{i}\lambda^{c_{i}(1+\phi_{i}+\cdots+\phi_{i}^{q_{i}-1}%
)}t^{\phi_{i}^{q_{i}}}-t
\]
and for $\underline{\lambda},\,\mathbf{t}\in F^{(M)}$ define
\begin{equation}
f_{\underline{\lambda}}(\mathbf{t})=\sum_{i=1}^{M}f_{i,\lambda_{i}}(t_{i}).
\label{sureq}%
\end{equation}
\medskip

\emph{(a)} Assume that $M>q(cq+1)$ and $\left|  F\right|  >c(cq+1)^{q}$ where
$c=\max\{c_{1},\ldots,c_{M}\}$. Then there exist $\lambda_{1},\ldots
,\lambda_{M}\in F^{\ast}$ such that the map $f_{\underline{\lambda}}%
:F^{(M)}\rightarrow F$ is surjective. \medskip

\emph{(b)} Let $\overline{F}\subseteq F$ be a subfield such that
$(F:\overline{F})=2$. Assume that $M>2q(2cq+1)$ and that $\left|  F\right|
>c^{2}(2cq+1)^{2q}$. Then the conclusion of part \emph{(a)} holds for some
$\lambda_{1},\ldots,\lambda_{M}\in\overline{F}^{\ast}$.
\end{lemma}

\begin{proof}
We shall give the proof of part (a) first and afterwards indicate the
modifications necessary to deduce part (b). \medskip

It will clearly suffice to show that for some subset $J\subseteq
\{1,\ldots,M\}$ and some $\lambda_{j}\in F^{\ast}$ ($j\in J$) the map
$F^{(\left|  J\right|  )}\rightarrow F$ given by
\[
(t_{i})_{i\in J}\mapsto\sum_{i\in J}f_{i,\lambda_{i}}(t_{i})
\]
is surjective. When $\left|  J\right|  =1$, this is equivalent to showing that
it has zero kernel (as it is linear over the prime field). Letting $F_{i}$
denote the fixed field of $\phi_{i}$ in $F$, we consider three cases.\medskip

\textbf{Case 1.} Suppose that for some $i$ we have $\mu_{i}\notin\lbrack
F^{\ast},\phi_{i}^{q_{i}}]$. Then for $t\in F^{\ast}$ we have
\[
f_{i,1}(t)=\mu_{i}t^{\phi_{i}^{q_{i}}}-t=\left(  \mu_{i}-[t^{-1},\phi
_{i}^{q_{i}}]\right)  t^{\phi_{i}^{q_{i}}}\neq0,
\]
so the kernel of $f_{i,1}:F\rightarrow F$ is zero; we take $J=\{i\},$
$\lambda_{i}=1$.

Henceforth, we may assume that for each $i$ there exists $b_{i}\in F^{\ast}$
such that
\[
\mu_{i}=b_{i}b_{i}^{-\phi_{i}^{q_{i}}}.
\]

\textbf{Case 2.} Suppose that $|F_{i}|>cq+1$ for some $i$. Let $\lambda$ be a
generator for the cyclic group $F^{\ast}$. We claim that the map
$f_{i,\lambda}:F\rightarrow F$ has zero kernel. Indeed, suppose not. Then
there exists $t\in F^{\ast}$ such that
\[
\mu_{i}\lambda^{c_{i}(1+\phi_{i}+\cdots+\phi_{i}^{q_{i}-1})}t^{\phi_{i}%
^{q_{i}}}=t,
\]
so writing $B=[F^{\ast},\phi_{i}]$ we have
\[
\lambda^{c_{i}(1+\phi_{i}+\cdots+\phi_{i}^{q_{i}-1})}=b_{i}^{-1}t\cdot
(b_{i}^{-1}t)^{-\phi_{i}^{q_{i}}}\in B.
\]
As $\lambda^{c_{i}(1+\phi_{i}+\cdots+\phi_{i}^{q_{i}-1})}\equiv\lambda
^{c_{i}q_{i}}$ modulo $B$ it follows that the cyclic group $F^{\ast}/B$ has
order dividing $c_{i}q_{i}$. But $\left|  F^{\ast}\right|  /\left|  B\right|
=\left|  F_{i}^{\ast}\right|  >cq\geq c_{i}q_{i}$, a contradiction. Thus we
may take $J=\{i\},$ $\lambda_{i}=\lambda$ in this case.\medskip

\textbf{Case 3.} Suppose that $|F_{i}|\leq cq+1$ for all $i$. Then there are
at most $cq+1$ possibilities for the size of each $F_{i}$, and also at most
$q$ possibilities for each $q_{i}$. As $M>q(cq+1)$ there exist $i<j\in
\{1,\ldots,M\}$ such that $|F_{i}|=|F_{j}|$ and $q_{i}=q_{j}$. Say $i=1$ and
$j=2$. We claim that in this case, there exists $\lambda\in F^{\ast}$ such
that the map
\[
(t_{1},t_{2})\mapsto f_{1,\lambda}(t_{1})+f_{2,1}(t_{2})
\]
from $F^{(2)}$ to $F$ is surjective.

It will suffice to show that for a suitable $\lambda$,
\[
(t_{1},t_{2})\mapsto b\left(  \lambda^{c_{1}(1+\phi_{1}+\cdots+\phi_{1}%
^{q_{1}-1})}t_{1}^{\phi_{1}^{q_{1}}}-t_{1}\right)  -\left(  t_{2}^{\phi
_{2}^{q_{2}}}-t_{2}\right)  =g_{\lambda}(t_{1},t_{2}),
\]
say, is surjective, where $b=-b_{1}b_{2}^{-1}$: replace $t_{i}$ by $b_{i}%
^{-1}t_{i}$ and divide by $-b_{2}$.

Since $\phi_{1}$ and $\phi_{2}$ have the same fixed field, they generate the
same subgroup of $\mathrm{Aut}(F)$. So $\phi_{1}=\phi_{2}^{l}$ for some $l$.
Writing
\[
t^{\ast}=t+t^{\phi_{2}^{q_{2}}}+\cdots+t^{\phi_{2}^{(l-1)q_{2}}}%
\]
and recalling that $q_{1}=q_{2}$ we have
\begin{align*}
g_{\lambda}(t,(bt)^{\ast})  &  =b\left(  \lambda^{c_{1}(1+\phi_{1}+\cdots
+\phi_{1}^{q_{1}-1})}t^{\phi_{1}^{q_{1}}}-t\right)  -\left(  (bt)^{\phi
_{1}^{q_{1}}}-bt\right) \\
&  =\left(  b\lambda^{c_{1}(1+\phi_{1}+\cdots+\phi_{1}^{q_{1}-1})}-b^{\phi
_{1}^{q_{1}}}\right)  t^{\phi_{1}^{q_{1}}}.
\end{align*}
It follows that $g_{\lambda}$ is surjective unless the bracketed factor is zero.

Suppose that this factor is zero for every $\lambda\in F^{\ast}$. Taking
$\lambda=1$ gives $b^{-1}b^{\phi_{1}^{q_{1}}}=1$, which then implies that
$\lambda^{c_{1}(1+\phi_{1}+\cdots+\phi_{1}^{q_{1}-1})}=1$ for every
$\lambda\in F^{\ast}$. Then $\lambda^{c_{1}\phi_{1}^{q_{1}}}=\lambda^{c_{1}}$,
hence $\lambda^{c_{1}}$ is in the fixed field $E$ of $\phi_{1}^{q_{1}}$ for
every $\lambda$, and it follows that $\left|  F^{\ast}\right|  \leq
c_{1}\left|  E^{\ast}\right|  $. On the other hand, the degree of the
extension $E/F_{1}$ divides $q_{1}$; as $|F_{1}|\leq cq+1$ this implies
$|E|\leq(cq+1)^{q}$. Thus
\[
\left|  F\right|  \leq1+c_{1}((cq+1)^{q}-1)\leq c(cq+1)^{q},
\]
contradicting the hypothesis.

This completes the proof of part (a). For part (b) we need some additional
notation: \medskip

Let $\overline{F}_{i}=\overline{F}\cap F_{i}$ be the fixed field of $\phi_{i}$
in $\overline{F}$ and let $\bar{\phi}_{i}$ be the generator of
Gal($F/\overline{F}_{i}$). If $\theta$ is the nontrivial automorphism of
$F/\overline{F}$ then $\langle\bar{\phi}_{i}\rangle=\langle\phi_{i}%
,\theta\rangle$. \medskip

The proof proceeds on the same lines as Part (a) with the following modifications:

\textbf{Case 1} is the same.

In \textbf{Case 2} we assume $|\overline{F}_{i}|>2cq+1$ and take $\lambda
:=\mu^{1+\theta}\in\overline{F}^{\ast}$, where $\mu$ is a generator of
$F^{\ast}$. Suppose that the kernel of $f_{i,\lambda}$ is non-zero. Define
$B:=[F^{\ast},\bar{\phi}_{i}]$. By considering $\lambda^{c_{i}(1+\phi
_{i}+\cdots+\phi_{i}^{q_{i}-1})}$ modulo $B$ we deduce that the cyclic group
$F^{\ast}/B$ has order $|\overline{F}_{i}^{\ast}|$ dividing $2c_{i}q_{i}$, a contradiction.

In \textbf{Case 3} we assume that $|\overline{F}_{i}|\leq2cq+1$ for all $i$.
As $F_{i}/\overline{F}_{i}$ has degree 1 or 2 this leaves at most $2(2cq+1)$
possibilities for $F_{i}$ and hence there is a pair $1\leq i<j\leq M$ such
that $F_{i}=F_{j}$ and $q_{i}=q_{j}$. We proceed as in Part (a) except that at
the end we reach the conclusion that for each $\lambda\in\overline{F}^{\ast}$
the element $\lambda^{c_{1}}$ is in the fixed field $\overline{E}%
\subseteq\overline{F}_{1}$ of $\phi_{1}^{q_{1}}\mid_{\overline{F}_{1}}$. This
gives $|\overline{F}^{\ast}|\leq c_{1}|\overline{E}^{\ast}|$ and as before
$|\overline{E}|\leq|\overline{F}_{1}|^{q_{1}}$. Therefore
\[
|F|=|\overline{F}|^{2}\leq(c|\overline{F}_{1}|^{q})^{2}\leq c^{2}%
(2cq+1)^{2q},
\]
contradicting the hypothesis of Part (b).
\end{proof}

\section{Orbital subgroups\label{orbitalsec}}

We can now give the\medskip

\textbf{Proof of Proposition \ref{orb}. }Recall that $S$ is a quasisimple
group of Lie type over a field $F$ of size at least $K=K(q)$, and that the
equivalence class $\omega\subseteq\Sigma_{+}$ is spanned by some orbit
$\overline{\omega}$ of positive roots from the untwisted system $\Sigma$. Thus
$\omega$ is the positive part of some (possibly orthogonally decomposable)
root system of dimension at most 3. In fact the type of $\omega$ is one of
$A_{1},A_{1}\times A_{1},A_{1}\times A_{1}\times A_{1},A_{2},B_{2}$ or $G_{2}%
$. Therefore $\omega$ has a height function with respect to its fundamental
roots. Let $\omega(i)$ be the set of (untwisted) roots of height at least $i $
in $\omega$.

The chief obstacle to applying Lemma \ref{sur} to $O=O(\omega)$ is that $O$
may not be abelian. However it is a nilpotent $p$-group and our strategy is to
find a filtration $O=O(1)\geq O(2)\geq\ldots\geq1$ of normal subgroups such
that each quotient $O(i)/O(i+1)$ is abelian and is even a vector space over
$F_{0}$ or $F$. This filtration is in fact provided by the sets $\omega(i)$ above.

For clarity we shall divide the proof in two parts, dealing first with the
case when $S$ is of untwisted type and second with the case when $S$ has
twisted type.

\subsection*{The untwisted case}

Recall that in this case we have $O(\omega)=W_{\omega}=\prod_{w\in\omega}%
X_{w}$. Define
\[
W^{i}=\prod_{v\in\omega(i)}X_{v}.
\]
This provides a natural filtration of subgroups
\begin{equation}
O=W_{\omega}=W^{1}\geq W^{2}\geq W^{3}\geq W^{4}=\{0\} \label{filtr}%
\end{equation}
of length at most 3. The factors $W^{i}/W^{i+1}$ are abelian and are modules
for the group $\mathcal{D}\Phi\Gamma$.

The automorphisms $\gamma_{i}^{q_{i}}$ may involve a graph automorphism and
thus may not stabilize the root subgroups. \medskip

We shall show that, provided $L_{0}$ and $|F|$ are sufficiently large, given
$i\in\{1,2,3\}$ and $\gamma_{1},\ldots,\gamma_{L_{0}}\in\mathcal{D}\Phi\Gamma
$, we can find elements $h_{j}\in H$ such that
\begin{equation}
W^{i}/W^{i+1}=\prod_{j=1}^{L_{0}}[W^{i}/W^{i+1},\beta_{j}],
\label{beta-module}%
\end{equation}
where $\beta_{j}=(\overline{h_{j}}\gamma_{j})^{e_{j}q_{j}}$ for some $e_{i} \in \{1,2,3\}$ to be chosen below. 
Then Lemma
\ref{trick} implies that the same will hold when each exponent $e_{j}q_{j}$ is
replaced by $q_{j}$, and Proposition \ref{orb} will follow, in view of
(\ref{filtr}), if we take $L\geq3L_{0}$.\bigskip

To establish (\ref{beta-module}), fix a root $w\in\omega$ and define $e_j=e_j (w)$ 
to be the size of the orbit of $w$ under the graph component of $\gamma_{j}$. Then 
$e_j \in \{1,2,3\}$ and
$\gamma_{j}^{e_{j}}$ stabilizes the root subgroup $X_w$. Set
\[
\alpha_{j}=(\overline{h_{j}}\gamma_{j})^{e_{j}}=\overline{h_{j}}^{1+\gamma
_{j}^{-1}+\cdots+\gamma_{j}^{1-e_{j}}}\gamma_{j}^{e_{j}},
\]
where the $h_{j}\in H$ remain to be determined. We shall show that for
a suitable $L_{1}$ it is possible to choose elements $h_{j}\in H$ so that
\begin{equation}
X_{w}\subseteq\prod_{j=1}^{L_{1}}[X_{w},\alpha_{j}^{q_{j}}]. \label{rootgp}%
\end{equation}
Since $\omega$ contains at most 6 roots, this will give (\ref{beta-module})
with $L_{0}=6L_{1}$ (on relabelling $\gamma_{j}$ and $e_{j}$, 
( $j=1,\ldots,L_{0})$ as $\gamma_{l}(w)$ and $e_l(w)$ with $w\in\omega,\,l=1,\ldots
,L_{1}$). \medskip

If $e_{j}=1$ we simply choose $h_{j}=h_{w}(\lambda_{j} )$ which acts on
$X_{w}$ as $t\mapsto\lambda_{j}^{2}t$.\medskip

If $e_{j}=2$ then the graph component $\tau$ of $\gamma_{j}^{-1}$ sends $w$ to
another root $v\in\omega$. Let
\[
h_{j}=h_{w}(\lambda_{j}^{2})h_{v}(\lambda_{j}^{-\langle v,w\rangle}).
\]
Then $h_{j}$ acts trivially on $X_{v}$, and on $X_{w}$ it acts as
$t\mapsto\lambda_{j}^{4-\langle\tau,v\rangle\langle v,w\rangle}t$. Notice that
$\langle w,v\rangle\langle v,w\rangle\in\{0,1,2,3\}$. Also $h_{j}^{\gamma
_{j}^{-1}}$ acts trivially on $X_{w}$. It follows that $h_{j}h_{j}^{\gamma
_{j}^{-1}}$ acts on $X_{w}$ as multiplication by $\lambda_{j}^{c_{j}}$, where
$c_{j}\in\{1,2,3,4\}$. \medskip

If $e_{j}=3$, set $h_{j}=h_{w}(\lambda_{j})$. The roots $w,w^{\tau}%
,w^{\tau^{2}}$ are pairwise orthogonal hence $h_{j}$ acts trivially on
$X_{w^{\tau}}$ and $X_{w^{\tau^{2}}}$, i.e. $h_{j}^{\tau}$ and $h_{j}%
^{\tau^{2}}$ act trivially on $X_{w}$. Therefore $h_{j}^{1+\gamma_{j}%
^{-1}+\gamma_{j}^{-2}}$ acts on $X_{w}$ as multiplication by $\lambda_{j}^{2}%
$. \medskip

Suppose $\gamma_{j}^{e_{j}} \in \mathcal{D} \Phi \Gamma$ acts on $X_w(t)$ as $t \mapsto \nu_j t^{\phi_j}$ with 
$\nu_{j}\in F^*$
and $\phi_{j}\in \mathrm{Aut}(F)$. Then an easy computation shows that
$\alpha_{j}^{q_{j}}=\left(  \overline{h_{j}}^{1+\cdots+\gamma_{j}^{1-e_{j}}%
}\gamma_{j}^{e_{j}}\right)  ^{q_{j}}$ acts on $X_{w}(t)$ as
\begin{equation}
t\mapsto\mu_{j}\lambda_{j}^{c_{j}\phi_{j} (1+\phi_{j}+\cdots+\phi_{j}%
^{q_{j}-1})}t^{\phi_{j}^{q_{j}}}, \label{action}%
\end{equation}
where $c_{j}\in\{1,2,3,4\}$ and $\mu_{j}$ is a constant which depends on
$\nu_{j},\phi_{j}$. In each case, therefore, (\ref{rootgp}) follows from Lemma
\ref{sur}, with $\lambda_{j}^{\phi_{j}}$ in place of $\lambda_{j}$, as long as
we take $L_{1}>q(4q+1)$ and $K>4(4q+1)^{q}$.

\subsection*{The twisted case}

Suppose now that $S=S^{\ast}$ is twisted with root system $\Sigma^{\ast}$
coming from the corresponding untwisted root sytem $\Sigma$. The complication
here is that the root subgroups are not necessarily 1-parameter. They are
parametrized by the same equivalence classes $\omega$ of roots under the
equivalence relation $\sim$ of $\Sigma$ considered above. Recall the
definition the root subgroups $Y_{\omega}$ and their description in Section
\ref{resume}. \medskip

Notice that $Y_{\omega}$ can be considered as a subgroup of $W_{\omega}$ in
the untwisted group above defined for $\Sigma$: in fact it is the subgroup of
$W_{\omega}$ fixed by the Steinberg automorphism $\sigma$ and inherits a
filtration $Y_{\omega}=Y^{1}\geq Y^{2}\geq Y^{3}\geq Y^{4}=\{0\}$ from $W$. We
shall therefore take $L>3L_{2}$ and show that for large enough $L_{2}$ it is
always possible to choose $\beta_{j}=(\overline{h_{j}}\gamma_{j})^{q_{j}}$
with $h_{j}\in H$ so that
\begin{equation}
Y^{i}/Y^{i+1}=\prod_{j=1}^{L_{2}}[Y^{i}/Y^{i+1},\beta_{j}]. \label{Y-module}%
\end{equation}

It is straightforward to adapt the strategy from the previous section since
the parametrization of $Y_{\omega}$ is coming ready from the ambient untwisted
group. However we keep in mind that we don't have all the diagonal elements at
our disposal: only those fixed by $\sigma$.

\medskip

In all cases $Y^{i}/Y^{i+1}$ is a one-parameter group, parametrized by either
$F_{0}$ or $F$. Also now $\Gamma=1$, and we shall apply the argument from the
previous subsection with each $e_{j}=1$. \medskip

\emph{The case of} $F_{0}$: When $\omega=A_{1}$ and $\,i=1$, or $\omega=A_{2}$
and $\,i=2$, the set $\omega(i)\setminus\omega(i+1)$ contains a single root
$v$. Then
\[
Y^{i}=X_{v}(aF_{0})\cdot Y^{i+1},
\]
where $a=1$ for type $A_{1}$, and $a\in F$ is any solution of $a+a^{\phi}=0$
for type $A_{2}$.

In this case, we set $h_{j}=h_{v}(\lambda_{j})\in H$ where the $\lambda_{j}\in
F_{0}$ are chosen so that the map $f_{\underline{\lambda}}$ defined by
(\ref{sureq}), with $c_{i}=2$ for each $i$, is surjective. This is possible by
Lemma \ref{sur} (applied to $F_{0}$) provided $|F_{0}|>2(2q+1)^{q}$ and
$L_{2}>q(2q+1)$ \ Then (\ref{Y-module}) follows as in the preceding
subsection.\medskip

\emph{The case of} $F$: In all other cases except for $\omega= (A_{1})^{3}$,
the set $\omega(i)\setminus\omega(i+1)=\{v,w\}$ consists of a pair of roots,
swapped by $\tau$. Say $v$ is the longer root. Then
\[
Y^{i}=\left\{  X_{v}(\xi)X_{w}(\xi^{\phi})\mid\xi\in F\right\}  \cdot
Y^{i+1}.
\]

Again we argue as in the preceding subsection (with $e_{j}=1$), but this time
we set
\[
h_{j}=h_{v}(\rho_{j})h_{w}(\rho_{j}^{\phi}),
\]
for suitably chosen $\rho_{j}\in F$. Then $h_{j}$ acts on $Y^{i}/Y^{i+1}$ via
$\xi\mapsto\rho_{j}^{2+\langle v,w\rangle\phi}\xi$. Notice that $|\langle
v,w\rangle|\in\{0,1,2,3\}$, and that if $\langle v,w\rangle\not =0$ then
$[|\langle v,w\rangle|]\phi^{2}$ is the indentity automorphism of $F$.
Consequently
\[
t^{4-(\langle v,w\rangle\phi)^{2}}=t^{4-|\langle v,w\rangle|} \quad\forall\ t
\in F.
\]

Taking $\rho_{j}=\lambda_{j}^{2-\langle v,w\rangle\phi}$ makes $h_{j}$ act on
$Y^{i}/Y^{i+1}$ via
\[
\xi\mapsto\lambda_{j}^{4-|\langle v,w\rangle|}\cdot\xi.
\]
So again we can apply Lemma \ref{sur} with $c_{i}=c=4-|\langle v,w\rangle|$
for each $i$, provided we assume that $L_{2}>q(4q+1)$ and $|F|>4(4q+1)^{q}$;
and (\ref{Y-module}) follows as above.

Finally, when $\omega=\{v,v^{\tau},v^{\tau^{2}}\}$ is of type $(A_{1})^{3}$ we
use
\[
h_{j}=h_{v}(\lambda_{j})h_{v^{\tau}}(\lambda_{j}^{\phi_{j}})h_{v^{\tau^{2}}%
}(\lambda_{j}^{(\phi_{j})^{2}}),
\]
($\lambda_{j} \in F$). This acts on $Y_{\omega}(t)$ as $t \mapsto\lambda
_{j}^{2} t$, and we apply Lemma \ref{sur} with all $c_{i}=2$.

\section{The unitriangular group\label{unitri}}

Here we establish Propositions \ref{SL} and \ref{manyW}. We begin with the
latter which concerns the group $V=V_{s+1}$ of unitriangular matrices in
$S=\mathrm{SL}_{s+1}(F)$ or $\mathrm{SU}_{s+1}(F)$ that differ from the
identity only in the first row and last column; here $s\geq5$ and $\left|
F\right|  >K=K(q)$.

We shall consider only the unwisted case $S=\mathrm{SL}_{s+1}(F)$. If
$S=\mathrm{SU}_{s+1}(F)$ the proof proceeds in exactly the same way by
considering the fixed points of $\sigma$ on the groups $V^{i}$ (and there is
no need to square $(\overline{h}_{j}\gamma_{j})^{q_{j}}$).

Now, the group $V$ has a filtration
\[
V=V^{1}>V^{2}>V^{3}>1
\]
where $V^{2}=\{g\in V\mid g_{12}=g_{s,s+1}=0\}$ and $V^{3}=\{g\in V\mid
g_{1j}=g_{i,s+1}=0\,$\ for $1<j<s,\,2<i<s+1\}$. Write
\[
W_{i}=\left\langle \mathbf{1}+e_{1,i+1}F,\,\mathbf{1}+e_{s+1-i,s+1}%
F\right\rangle
\]
where $e_{ij}$ denotes the matrix with $1$ in the $(i,j)$ entry and $0$
elsewhere. Each $W_{i}$ is an orbital subgroup of $S$ and we have
\[
V=W_{1}\cdot V^{2},\,V^{2}=W_{2}W_{3}\ldots W_{s-2}\cdot V^{3},\,\,\,V^{3}%
=W_{s-1}.
\]

Now let $\gamma_{1},\gamma_{2},\ldots$ be automorphisms of $S$ lying in
$\mathcal{D}\Phi\Gamma$ and let $q_{1},q_{2},\ldots$ be divisors of $q$. We
have to find elements $h_{1},h_{2},\ldots\in H$ such that
\[
V=\prod_{i=1}^{L_{1}}[V,(\overline{h_{i}}\gamma_{i})^{q_{i}}].
\]
Let $L=L(q)$ be the integer given in Proposition \ref{orb}, and take
$L_{1}=2L+L_{3}$. Applying that proposition to $W_{1}$ and to $W_{s-1},$ and
relabelling the $\gamma_{i}$ and $q_{i}$, we are reduced to showing that there
exist $h_{1},\ldots,h_{L_{3}}\in H$ such that
\begin{equation}
V^{2}/V^{3}=\prod_{i=1}^{L_{3}}[V^{2}/V^{3},(\overline{h_{i}}\gamma
_{i})^{q_{i}}]; \label{V2/V3}%
\end{equation}
note that $V^{2}/V^{3}\cong F^{(2(s-3))}$ is a module for $\mathcal{D}%
\Phi\Gamma$.

Now for any $\lambda_{i}\in F^{\ast}$ we may choose $h_{i}\in H$ so that the
diagonal component of $\overline{h_{i}}\gamma_{i}$ is of the form
$\mathrm{diag}(\lambda_{i}^{-1},\ast,1,\ldots,1,\ast,\lambda_{i})$; this acts
on $V^{2}/V^{3}$ as multiplication by $\lambda_{i}$. Let $\phi_{i}$ denote the
field component of $\gamma_{i}$. Provided $L_{3}>2q(2q+1)$, Lemma \ref{sur}
gives elements $\lambda_{1},\ldots,\lambda_{L_{3}}\in F^{\ast}$ such that the
map $f : F^{(L_{3})}\rightarrow F$ defined by
\[
(t_{1},t_{2},\ldots,t_{L_{3}})\mapsto\sum_{i=1}^{L_{3}}(\lambda_{i}^{\phi
_{i}(1+\phi_{i}+\cdots+\phi_{i}^{2q_{i}-1})}t_{i}^{\phi_{i}^{2q_{i}}}-t_{i})
\]
is surjective. Formula (\ref{action}), from the previous section, shows that
the action of $\sum_{i} ((\overline{h_{i}}\gamma_{i})^{2q_{i}}-1)$ on each
root of $V^{2}/V^{3}$ is in fact given by $f$. Therefore
\[
V^{2}/V^{3}=\prod_{i=1}^{L_{3}}[V^{2}/V^{3},(\overline{h_{i}}\gamma
_{i})^{2q_{i}}],
\]
and (\ref{V2/V3}) follows by Lemma \ref{trick}.

This completes the proof of Proposition \ref{manyW}, with $L_{1}%
(q)=2L(q)+2q(2q+1)+1$.\bigskip

It remains to prove Proposition \ref{SL}. Let $S=\mathrm{SL}_{r+1}(F)$, where
$\left|  F\right|  >K$ and $r\geq3$, and put $M_{2}=4L(q)+1$. We are given
automorphisms $\gamma_{1},\ldots,\gamma_{M_{2}}\in\mathcal{D}\Phi\Gamma$ and
divisors $q_{1},\ldots,q_{M_{2}}$ of $q,$ and have to find automorphisms
$\eta_{1},\ldots,\eta_{M_{2}}\in\mathcal{D}$ and elements $u_{1}%
,\ldots,u_{M_{2}}\in U$ such that
\begin{equation}
U\subseteq\prod_{i=1}^{M_{2}}[U,(\overline{u_{i}}\eta_{i}\gamma_{i})^{q_{i}}]
\label{U}%
\end{equation}
(here $U$ is the full upper unitriangular group in $S$). Since we are allowed
to adjust $\gamma_{i}$ by any element of $\mathcal{D}$, we may without loss of
generality assume from now on that each $\gamma_{i}\in\Phi\Gamma$.

A matrix $g\in U$ will be called \emph{proper} if $g_{i,i+1}\not =0$ for
$1\leq i\leq r$. Let $U_{k}$ be the product of all positive root groups of
height $\geq k$ (so $u\in U_{k}$ precisely if $u_{ij}=0$ for $0<j-i<k$). It is
well known that $U_{1}>U_{2}>\ldots$ is the lower central series of $U$ and
that for each $k\leq r-1$ and each proper matrix $g$ the map $x\mapsto\lbrack
x,g]$ induces a surjective linear map of $\mathbb{F}_{p}$-vector spaces
\begin{equation}
\lbrack-,g]:U_{k}/U_{k+1}\rightarrow U_{k+1}/U_{k+2}. \label{linmap}%
\end{equation}
We call the section $U_{k}/U_{k+1}$ the $k$\emph{-th layer of U}. By a slight
abuse of notation, we shall identify $\mathcal{D}$ with the group of diagonal
matrices in $\mathrm{GL}_{r+1}(F)$ modulo scalars.

\begin{lemma}
\label{g0} Let $q_{0}\in\mathbb{N}$ and suppose that $|F|>(q_{0}+1)^{q_{0}}$.
Then for any $\gamma\in\Phi\Gamma$ we can find $u\in U$ and $\eta
\in\mathcal{D}$ \ such that the matrix
\[
g=(u\eta\gamma)^{q_{0}}(\eta\gamma)^{-q_{0}}%
\]
is proper.
\end{lemma}

\begin{proof}
Let $\gamma\in\Phi\Gamma$ act on the root subgroup $X_{r}(t)$ as
$X_{r}(t)^{\gamma}=X_{r^{\tau}}(t^{\psi})$ (here $\psi$ is a field
automorphism of $F,$ and $\tau$ may be $1$).\medskip

\textbf{Case 1}: When $\psi^{q_{0}}\not =1$. Then we can find $\lambda\in
F^{\ast}$ such that $\lambda\not =\lambda^{\psi^{q_{0}}}$, and so
\[
\mu=\lambda^{1+\psi^{-1}+\cdots+\psi^{1-q_{0}}}\not =1.
\]
Let $h(\lambda)$ be the diagonal automorphism $\mathrm{diag}(1,\lambda
^{-1},\lambda^{-2},\ldots)$, acting on each fundamental root group $X_{r}$
($r\in\Pi$) by $t\mapsto\lambda t$. Then $h(\lambda)$ commutes with $\tau$ and
we have $h(\mu)\gamma^{q_{0}}=(h(\lambda)\gamma)^{q_{0}}$.

Now let $v=\prod_{r\in\Pi}X_{r}(1)$ be the unitriangular matrix with $1$s just
off the diagonal. Then $v$ is centralized by $\gamma$ modulo $U_{2}$ and
$[v,h(\mu)^{-1}]$ is proper. Putting
\[
u=[v,h(\lambda)^{-1}],\,\,\eta=h(\lambda)
\]
we have, modulo $U_{2}$,
\[
(u\eta\gamma)^{q_{0}}\equiv(h(\lambda)^{v}\gamma^{v})^{q_{0}}\equiv\left(
(h(\lambda)\gamma)^{q_{0}}\right)  ^{v}\equiv(h(\mu)\gamma^{q_{0}})^{v}%
\equiv\lbrack v,h(\mu)^{-1}](\eta\gamma)^{q_{0}},
\]
so $g\equiv\lbrack v,h(\mu)^{-1}]$ is proper as required. \medskip

\textbf{Case 2}: When $\psi^{q_{0}}=1$. Let $F_{1}$ be the fixed field of
$\psi$. Then $[F:F_{1}]\leq q_{0}$, so if $|F|>(q_{0}+1)^{q_{0}}$ it follows
that $|F_{1}|>q_{0}+1$. Therefore we can choose $\lambda\in F_{1}$ such that
\[
\lambda^{q_{0}}=\lambda^{1+\psi^{-1}+\cdots+\psi^{1-q_{0}}}\not =1,
\]
and the rest of the proof is as in Case 1.
\end{proof}

\bigskip

To establish (\ref{U}), we begin by showing that we can obtain the slightly
smaller group $U_{3}$ as a product of $2L+1$ sets $[U_{3},(\overline{u_{i}%
}\eta_{i}\gamma_{i})^{q_{i}}]$.

\begin{lemma}
\label{U3}Let $\gamma_{0},\gamma_{1},\ldots,\gamma_{2L}\in\Phi\Gamma$ and let
$q_{0},\ldots,q_{2L}$ be divisors of $q$. Then there exist $\eta_{i}%
\in\mathcal{D}\ \ (i=0,\ldots,2L$) and $u\in U$ such that
\[
U_{3}=\prod_{i=1}^{2L}[U_{3},(\eta_{i}\gamma_{i})^{q_{i}}]\cdot\lbrack
U_{3},(\overline{u}\eta_{0}\gamma_{0})^{q_{0}}].
\]

\end{lemma}

\noindent(So here we have $u_{i}=1$ for $i=1,\ldots,2L$ and $u_{0}=u$.)

\begin{proof}
First, note that if $\alpha\in\mathrm{Aut}(U)$ and $x,y\in U_{k}$ then
\begin{equation}
\lbrack xy,\alpha]=[x,\alpha][x,\alpha,y][y,\alpha]\equiv\lbrack
x,\alpha][y,\alpha]\,\,(\operatorname{mod}U_{2k}). \label{xy}%
\end{equation}

Say $\gamma_{i}=\phi_{i}$ or $\phi_{i}\tau$ where $\phi_{i}\in\Phi$. Then, by
a double application of Lemma \ref{sur}, we can find $\lambda_{i}\in F^{\ast
},\ i=1,2,\ldots,2L$, such that \emph{both} of the maps $f_{+},\,f_{-}%
:F^{(2L)}\rightarrow F$ defined by
\begin{align*}
f_{+}(\mathbf{t})  &  =\sum_{i=1}^{2L}\lambda_{i}^{\phi_{i}+\phi_{i}%
^{2}+\cdots+\phi_{i}^{2q_{i}}}t_{i}^{\phi_{i}^{2q_{i}}}-t_{i},\\
f_{-}(\mathbf{t})  &  =\sum_{i=1}^{2L}\lambda_{i}^{-(\phi_{i}+\phi_{i}%
^{2}+\cdots+\phi_{i}^{2q_{i}})}t_{i}^{\phi_{i}^{2q_{i}}}-t_{i}%
\end{align*}
are surjective. Indeed, it suffices to ensure that each of the maps
\[
\mathbf{t}\mapsto\sum_{i=1}^{L}\lambda_{i}^{\phi_{i}+\phi_{i}^{2}+\cdots
+\phi_{i}^{2q_{i}}}t_{i}^{\phi_{i}^{2q_{i}}}-t_{i}%
\]
and
\[
\mathbf{t}\mapsto\sum_{i=L+1}^{2L}\lambda_{i}^{-(\phi_{i}+\phi_{i}^{2}%
+\cdots+\phi_{i}^{2q_{i}})}t_{i}^{\phi_{i}^{2q_{i}}}-t_{i}%
\]
is surjective.

We now take
\[
\eta_{i}=\left\{
\begin{array}
[c]{ccc}%
\mathrm{diag}(\lambda_{i},1,\lambda_{i},1,\ldots,\lambda_{i},1) &  & (s\text{
odd})\\
&  & \\
\mathrm{diag}(\ldots,1,\lambda_{i},1,\lambda_{i},\underline{1},\lambda
_{i}^{-1},1,\lambda_{i}^{-1},1,\ldots) &  & (s\text{ even})
\end{array}
\right.  \text{,}%
\]
where in the even rank case the underlined unit $\underline{1}$ has the
central position $1+\frac{s}{2}$ on the diagonal of $\mathrm{SL}_{s+1}$.
\medskip

It is easy to see that if $w$ is a root of \emph{odd} height, then either
$X_{w}(t)^{\eta_{i}}=X_{w}(\lambda_{i}t)$ for all $i$, or else $X_{w}%
(t)^{\eta_{i}}=X_{w}(\lambda_{i}^{-1}t)$ for all $i$. Moreover, it follows
from the definition that $\eta_{i}^{\tau}=\eta_{i}$. Then the surjectivity of
$f_{+}$ and $f_{-}$ together imply that $X_{w}=\prod_{i=1}^{2L}[X_{w}%
,(\eta_{i}\gamma_{i})^{2q_{i}}]$ for all roots $w$ of odd height. (Note that
$\gamma_{i}^{2q_{i}}$ stabilizes every root of $\Sigma$). Hence when $k\geq3$
is \emph{odd} we have
\[
U_{k}/U_{k+1}=\prod_{i=1}^{2L}[U_{k}/U_{k+1},(\eta_{i}\gamma_{i})^{2q_{i}}].
\]
It follows from Lemma \ref{trick} that the product $\prod_{i=1}^{2L}%
[U_{3},(\eta_{i}\gamma_{i})^{q_{i}}]$ covers each odd layer of $U_{3}$.

To deal with the even layers we use the map (\ref{linmap}). Put $\beta
_{i}=(\eta_{i}\gamma_{i})^{q_{i}}\in\mathcal{D}\Phi\Gamma$ for each $i$. Now
take $b\in U_{3}$ and let $k\geq3$ be odd. Suppose that we have already found
$x_{i},y\in U_{3}$ such that
\[
b\equiv\prod_{i=1}^{2L}[x_{i},\beta_{i}]\cdot\lbrack y,\overline{g}\beta
_{0}]\text{ }(\operatorname{mod}U_{k}),
\]
where $g=(u\eta_{0}\gamma_{0})^{q_{0}}(\eta_{0}\gamma_{0})^{-q_{0}}$ is the
proper matrix provided by Lemma \ref{g0}.

We claim that there exist $x_{1}^{\prime},x_{2}^{\prime},\ldots,x_{2L}%
^{\prime},\,$ $y^{\prime}\in U_{k}$ such that
\begin{equation}
b\equiv\prod_{i=1}^{2L}[x_{i}x_{i}^{\prime},\beta_{i}]\cdot\lbrack yy^{\prime
},\overline{g}\beta_{0}]\text{ }(\operatorname{mod}U_{k+2}). \label{xy-equn}%
\end{equation}
By (\ref{xy}) this is equivalent to
\begin{equation}
\prod_{i=1}^{2L}[x_{i}^{\prime},\bar{\beta}_{i}]\cdot\lbrack y^{\prime
},\overline{g}\beta_{0}]\equiv b^{\prime}\,\,(\operatorname{mod}U_{k+2})
\label{x'y'-equn}%
\end{equation}
where $b^{\prime}=b\cdot\left(  \prod[x_{i},\beta_{i}]\cdot\lbrack
y,\overline{g}\beta_{0}]\right)  ^{-1}\in U_{k}$. Also
\[
\lbrack y^{\prime},\overline{g}\beta_{0}]=[y^{\prime},\beta_{0}]\cdot
\lbrack{y^{\prime}},g]^{\beta_{0}}.
\]

Let
\[
V_{1}=\left(  U_{k+2}\prod_{\mathrm{ht}(w)=k}X_{w}\right)  /U_{k+2},\qquad
V_{2}=\left(  U_{k+2}\prod_{\mathrm{ht}(w)=k+1}X_{w}\right)  /U_{k+2}.
\]
We identify $U_{k}/U_{k+2}$ with $V_{1}\oplus V_{2}$. The elementary abelian
$p$-group $V_{1}$ corresponds to the (odd) $k$th layer of $U$, while $V_{2}$
is the (even) ($k+1$)th layer.

Now, on the one hand the map $y^{\prime}\mapsto\lbrack y^{\prime},g]U_{k+2}$
($y^{\prime}\in V_{1}$) is a surjective linear map from $V_{1}$ onto $V_{2}$.
On the other hand, the argument above shows that the map $\mathbf{x}^{\prime
}\mapsto\prod_{i=1}^{2L}[x_{i}^{\prime},\beta_{i}]U_{k+2}$ \ ($\mathbf{x}%
^{\prime}\in V_{1}^{(2L)}$) is a surjective linear map from $V_{1}^{(2L)}$
onto $V_{1}$.

We can therefore solve the equation (\ref{x'y'-equn}) in $U_{k}/U_{k+2}$ in
the following way. Suppose $b^{\prime}=b_{1}+b_{2}$ with $b_{i}\in V_{i}$.
First choose $y^{\prime}\in V_{1}$ so that $[{y^{\prime}},g]=b_{2}^{\beta
_{0}^{-1}}$ and observe that $[y^{\prime},\beta_{0}]\in V_{1}$. Consider this
$y^{\prime}$ fixed and then choose appropriate $x_{i}^{\prime}\in V_{1}$ so
that
\[
\prod_{i=1}^{2L}[x_{i}^{\prime},\beta_{i}]=b_{1}-[y^{\prime},\beta_{0}].
\]

It follows by induction on the odd $k,$ starting with $k=3$, that we can solve
(\ref{xy-equn}) with $k=r+1$, and as $U_{r+1}=1$ this establishes the lemma.
\end{proof}

\bigskip

What remains to be done is to obtain the first two layers $U_{1}/U_{2}$ and
$U_{2}/U_{3}$. We shall need $2L$ more automorphisms $(\eta_{i}\gamma
_{i})^{q_{i}}$.

The set of roots of height $1$ is $\Pi$, and we denote by $\Xi$ the set of
roots of height $2$. We show first how to obtain $\prod_{w\in\Pi}X_{w}$.

For a choice of $\lambda_{i}\in F^{\ast},\,i=1,2,\ldots,L$, put
\[
\eta_{i}=\mathrm{diag}(1,\lambda_{i},\lambda_{i}^{2},\lambda_{i}^{3},\ldots).
\]
For each fundamental root $w$ we then have
\[
X_{w}(t)^{\eta_{i}}=X_{w}(\lambda_{i}t).
\]
In particular the restrictions of $\eta_{i}$ and $\eta_{i}^{\tau}$ to
$U_{1}/U_{2} $ are the same. As in the proof of Proposition \ref{manyW},
above, we may apply Lemma \ref{sur} to find $\lambda_{i}\in F^{\ast}$ such
that
\[
U_{1}/U_{2}\subseteq\prod_{i=1}^{L}[U_{1}/U_{2},(\eta_{i}\gamma_{i})^{2q_{i}%
}]\subseteq\prod_{i=1}^{L}[U_{1}/U_{2},(\eta_{i}\gamma_{i})^{q_{i}}].
\]
The analogous result for $\prod_{w\in\Xi}X_{w}$ is obtained similarly using
the diagonal automorphisms
\[
\eta_{i}=\mathrm{diag}(1,1,\lambda_{i},\lambda_{i},\lambda_{i}^{2},\lambda
_{i}^{2},\ldots).
\]

As $U=(\prod_{w\in\Pi}X_{w})\cdot(\prod_{w\in\Xi}X_{w})\cdot U_{3}$, the last
two observations together with Lemma \ref{U3} complete the proof of
Proposition \ref{SL}, with $M_{2}=4L+1$.

\section{The group $P$\label{Psection}}

Here we prove Propositions \ref{P1} and \ref{P2}. Let us recall the setup. $S$
is a quasisimple group of type%
\[
\mathcal{X}\in\{^{2}A_{r},\,B_{r},\,C_{r},\,D_{r},\,^{2}D_{r}\},
\]
with root system $\Sigma$ (twisted or untwisted). Here $r$ can be any integer
greater than $3$. There exist fundamental roots $\delta,\,\delta^{\prime}%
\in\Sigma$ (equal unless $\mathcal{X}=D_{r}$) such that the other fundamental
roots $\Pi^{\prime}=\Pi-\{\delta,\delta^{\prime}\}$ generate a root system
$\Sigma^{\prime}$ of type $A_{s}$, for the appropriate $s$: in types
$^{2}A_{r},B_{r},C_{r}$ and $^{2}D_{r}$ we take $\delta=\delta^{\prime}$ to be
the fundamental root of length distinct from the others; in type $D_{r}$,
$\{\delta,\delta^{\prime}\}$ is the pair of fundamental roots swapped by the
symmetry $\tau$ of $D_{r}$.

If $S$ is untwisted, set%
\[
P=\prod_{w\in\Sigma_{+}\backslash\Sigma_{+}^{\prime}}X_{w}.
\]
If $S$ is twisted, set%
\[
P=\prod_{\omega^{\ast}\in\Sigma_{+}\backslash\Sigma_{+}^{\prime}}Y_{\omega}.
\]

\noindent\textbf{Proposition \ref{P1} }\emph{Assume that } $\left|  F\right|
>K$\emph{ and that }$S$\emph{ is of type } $B_{r},\,C_{r},\,D_{r}$\emph{ or
}$\,^{2}D_{r}$\emph{. There is a constant }$N_{1}=N_{1}(q)$\emph{ such that if
}$\gamma_{1},\ldots,\gamma_{N_{1}}$\emph{ are automorphisms of }$S$ \emph{
lying in }$\mathcal{D}\Phi\Gamma$\emph{ and }$q_{1},\ldots,q_{N_{1}}$
\emph{are divisors of }$q$\emph{ then there exist elements } $h_{1}%
,\ldots,h_{N_{1}}\in H$\emph{ such that}
\[
P\subseteq\prod_{i=1}^{N_{1}}[P,(\overline{h_{i}}\gamma_{i})^{q_{i}}].
\]

\bigskip

We consider first the \textbf{untwisted} case, where $S$ has type $B_{r}$,
$C_{r}$ or $D_{r}$. By inspection of the root systems we see that every
positive root $w\in\Sigma_{+}$ can be written as $w=e\delta+w_{1}+e^{\prime
}\delta^{\prime}+w_{2}$ with some $e,e^{\prime}\in\{0,1\}$ and $w_{1},w_{2}%
\in\Sigma_{+}^{\prime}\cup\{0\}$. In the last expression we include the
possibility $\delta=\delta^{\prime}$.

For $w=(e\delta+w_{1})+(e^{\prime}\delta^{\prime}+w_{2})$ as above we set
$t(w):=e+e^{\prime}$.

For $i=1,2$ let $P(i)$ be the product of roots subgroups $X_{w}$ with
$t(w)\geq i$ in any order; this is in fact a normal subgroup of $P$. Then
$P=P(1)$, and both $P(1)/P(2)$ and $P(2)$ are abelian. Each of $P(1)$ and
$P(2)$ is a product of orbital subgroups and so invariant under $\mathcal{D}%
\Phi\Gamma$.

Recall Lemma \ref{groupD}. The type of $S$ is here different from $A_{n}$ and
$^{2}A_{n}$ and therefore there are characters $\chi_{i}:\Sigma\rightarrow
F^{\ast}$ \ ($i=1,\ldots,4$) \ of the root lattice such that (i)%
\[
\mathcal{D}=\bigcup_{i=1}^{4}h(\chi_{i})\overline{H}%
\]
and (ii)%
\[
\chi_{i}(w)=1\text{ }\forall w\in\Pi\setminus\Delta
\]
where $\Delta$ is a fixed set of fundamental roots of size at most 2: if
$X=D_{r}$ then $\Delta=\{\delta,\delta^{\prime}\}$, otherwise $\Delta$
consists of a single root at one end of the Dynkin diagram

In view of (i), we may assume that the diagonal component $d_{j}$ of each
$\gamma_{j}$ is one of the four $h(\chi_{i})$ above. Setting $N_{1}>4N_{2}$
and relabelling the $\gamma_{j}$ if necessary we may further suppose that%
\[
1\leq j\leq N_{2}\Longrightarrow d_{j}=h(\chi_{1})=h_{0},\,\text{say.}%
\]
Observe that for any root $w\in\Sigma_{+}$ the multiplicity of each of
$v\in\Delta$ in $w$ is $0$, $1$ or $2$, and the last case occurs only when
$\Delta$ has size $1$. Therefore $h_{0}$ can act in only $4$ possible
different ways on $X_{w}$ as $w$ ranges over $\Sigma_{+}$. We deduce that a
given $\gamma_{j}^{2q_{j}}$ can act in at most four different ways on the
various root subgroups. (The possible presence of a graph automorphism
component of $\gamma_{j}$ in case of $D_{r}$ requires additional attention.)
We make this more precise:

Given $\gamma_{j}$ for $j=1,2,\ldots,N_{2}$, each with diagonal component
$h_{0}$, for each $j$ there are elements $c_{j}(i)\in F^{\ast}$ $\ $%
($i=1,\ldots,4$) with the following property:\medskip

For each root $w\in\Sigma_{+}$ there exists $i=i(w)\in\{1,2,3,4\}$ such that
$\gamma_{j}^{2q_{j}}$ acts on $X_{w}$ as
\[
X_{w}(t)\mapsto X_{w}(c_{j}(i)t^{\phi_{j}^{2q_{j}}}).
\]
Here $\phi_{j}$ is the field automorphism component of $\gamma_{j}$. \medskip

For any $\lambda\in F^{\ast}$ let $\chi_{\lambda}$ denote the character of
$\Sigma$ which takes value $\lambda^{4}$ on $\{\delta,\delta^{\prime}\}$ and
is $1$ on all the other fundamental roots. The automorphism $h_{\lambda
}:=h(\chi_{\lambda})$ is a fourth power in $\mathcal{D}$ and therefore inner.
Observe that if $v\in\Sigma_{+}\backslash\Sigma_{+}^{\prime}$ then
\[
X_{v}(t)^{h_{\lambda}}=X_{v}(\lambda^{4\cdot t(v)}t),
\]
where $t(v)\in\{1,2\}$ was defined above.

Let $\phi_{j}$ be the field automorphism component of $\gamma_{j}$. The automorphisms
$(h_{\lambda_{j}}\gamma_{j})^{2q_{j}}$ stabilize each root subgroup $X_{w}\leq
P$ and act on $X_{w}(t)$ as
\[
t\mapsto c_{j}(i)\lambda_{j}^{4t(w)\cdot(\phi_{j}+\phi_{j}^{2}+\cdots+\phi
_{j}^{2q_{j}})}\cdot t^{\phi_{j}^{2q_{j}}},
\]
where $i=i(w)\in\{1,2,3,4\}$ as above, $t(w)=2$ if $X_{w}\leq P(2)$ and
$t(w)=1$ otherwise.

We set $N_{2}=N_{3}+N_{4}$. Let $\{1,2,\ldots,N_{2}\}=J_{3}\cup J_{4}$ where
$J_{3}$ and $J_{4}$ have sizes $N_{3}$ and $N_{4}$ respectively. set
$N_{3}=4M$ and let $J_{3}$ be a union of four subsets $J_{3}(i)$,
$i=1,\ldots,4$ each of size $M$.

Using Lemma $\ref{sur}$ with all $c_{i}=4$, provided $|F|$ is large compared
to $q$ and $M>2q(8q+1)$ we may find $\underline{\lambda}\in F^{(J_{3})}$ such
that the maps
\begin{equation}
\mathbf{t}\in F^{(J_{3}(i))}\mapsto\sum_{j\in J_{3}(i)}\left(  c_{j}%
(i)\lambda_{j}^{4(\phi_{j}+\phi_{j}^{2}+\cdots+\phi_{j}^{2q_{j}})}\cdot
t_{j}^{\phi_{j}^{2q_{j}}}-t_{j}\right)  \label{seventeen}%
\end{equation}
are surjective for each $i=1,2,3,4$.

This gives
\begin{equation}
P/P(2)\subseteq\prod_{j\in J_{3}}[P/P(2),(h_{\lambda_{j}}\gamma_{j})^{2q_{j}%
}]\subseteq\prod_{j\in J_{3}}[P/P(2),(h_{\lambda_{j}}\gamma_{j})^{q_{j}}]
\label{eighteen}%
\end{equation}
Similarly, another four-fold application of Lemma \ref{sur} with $c_{i}=8$
gives that for $N_{4}>4\times2q(16q+1)$ there exist $\lambda_{j}\in F^{\ast}$
for $j\in J_{4}$ such that the analogue of (\ref{seventeen}) holds and hence
\begin{equation}
P(2)\subseteq\prod_{j\in J_{4}}[P(2),(h_{\lambda_{j}}\gamma_{j})^{2q_{j}%
}]\subseteq\prod_{j\in J_{4}}[P(2),(h_{\lambda_{j}}\gamma_{j})^{q_{j}}]
\label{nineteen}%
\end{equation}
This concludes the proof in the untwisted case.

\bigskip

It remains to establish the \textbf{twisted} case, where $S$ has type
$^{2}D_{r}$. We will denote by $\Sigma^{0}$ the untwisted root system
corresponding to $\Sigma$.

The group $P=P^{\ast}$ inherits a filtration $P=P^{\ast}(1)>P^{\ast}(2)$ from
the associated untwisted root system $D_{r}$: each $P^{\ast}(i)$ is the fixed
points of $\sigma$ on the corresponding group $P(i)$ defined as above for the
untwisted version of $S$: the subgroup $P^{\ast}(i)$ is the product of all
root subgroups $Y_{\omega}$ with equivalence class $\omega$ consisting of
untwisted roots $w$ with $t(w)\leq i$. Recall that in type $^{2}D_{r}$ the
root subgroups $Y_{\omega}$ are all one-parameter.

The group $P^{*}(2)$ is the product of the root subgroups $Y_{\omega}$ defined
by a singleton $\omega=\{w\}$ where the positive root $w \in\Sigma^{0}$ is
fixed by $\tau$. Then $Y_{\omega}=\{X_{w}(t)\ \mid\ t\in F_{0}\}$.

On the other hand the group $P^{\ast}/P^{\ast}(2)$ is the product of
$Y_{\omega}P^{\ast}(2)$ where $\omega=\{u,v\}\subseteq\Sigma_{+}^{0}$ has type
$A_{1}\times A_{1}$ and $Y_{\omega}(t)=X_{u}(t)X_{v}(t^{\phi})$ is
parametrized by $t\in F$.

We now proceed as in the previous case:

By Lemma \ref{groupD} we may take $N_{1}>4N_{2}$ and may assume that the
automorphisms $\gamma_{j}$ all have the same diagonal component $h_{0}$ for
$j=1,2,\ldots,N_{2}$.

There are elements $c_{j}(i)\in F^{\ast}$, ($1\leq j\leq N_{2}$, $i=1,2,3,4$),
depending on $h_{0}$, such that the automorphism $\gamma_{j}^{q_{j}}$ acts on
a root element $Y_{\omega}(t)$ as $t\mapsto c_{j}(i)t^{\phi_{j}^{q_{j}}}$,
where $i=i(\omega)\in\{1,2,3,4\}$ depends only on $\omega$.

Let $\{a,b\}$ be the pair of fundamental roots in $\Sigma^{0}$ corresponding
to the short root $\delta\in\Sigma$. For $\lambda\in F^{\ast}$ let
$\chi_{\lambda}$ be the character of the untwisted root system $\Sigma^{0}$
defined by $\chi(a)=\lambda^{4}$, $\chi(b)=\lambda^{4\phi}$ and $\chi$ is 1 on
the rest of the fundamental roots of $\Sigma^{0}$.

Define $h_{\lambda}:=h(\chi)$. Then $h_{\lambda}$ is an inner diagonal
automorphism and is fixed by $\sigma$, therefore $h_{\lambda}\in\overline{H}$.
If $\phi_{j}$ is the field component of $\gamma_{j}$ then $(h_{\lambda_{j}%
}\gamma_{j})^{q_{j}}$ acts on $Y_{\omega}(t)$ as
\[
t\mapsto c_{j}(i)\lambda_{j}^{c_{\omega}(\phi_{j}+\cdots+\phi_{j}^{q_{j}}%
)}t^{\phi_{j}^{q_{j}}},
\]
where $c_{j}(i)$ and $i=i(\omega)$ are as above, and

\begin{itemize}
\item {$c_{\omega}=4+4\phi$ if $\omega$ has type $A_{1}$ (when $t$ ranges over
$F_{0}$),}

\item {$c_{\omega}=4$ if $\omega$ has type $A_{1}\times A_{1}$ (when $t$
ranges over $F$).}
\end{itemize}

We set $N_{2}=N_{3}+N_{4}$.

In the same way as in the previous case, provided $N_{3}$, $N_{4}$ and
$|F_{0}|$ are sufficiently large compared to $q$, it is possible to choose
$\lambda_{j}\in F^{\ast}$ for $j=1,2,\ldots,N_{3}$ and $\lambda_{j}\in
F_{0}^{\ast}$ for $N_{3}<j\leq N_{4}$ so that the appropriate equivalents of
(\ref{seventeen}) hold. This gives (\ref{eighteen}) and (\ref{nineteen}) and
concludes the proof in the case of type $^{2}D_{r}$.

\bigskip

\noindent\textbf{Proposition \ref{P2}} \emph{Assume that }$\left|  F\right|
>K$\emph{ and that }$S$ \emph{ is of type }$\,^{2}A_{r}$\emph{. There is a
constant }$N_{1}^{\prime}=N_{1}^{\prime}(q)$\emph{ such that if } $\gamma
_{1},\ldots,\gamma_{N_{1}^{\prime}}$\emph{ are automorphisms of }$S$ \emph{
lying in }$\mathcal{D}\Phi\Gamma$\emph{ and } $q_{1},\ldots,q_{N_{1}^{\prime}%
}$\emph{ are divisors of }$q$ \emph{ then there exist automorphisms }
$\eta_{1},\ldots,\eta_{N_{1}^{\prime}}\in\mathcal{D}$\emph{ such that}
\[
P\subseteq\prod_{i=1}^{N_{1}^{\prime}}[P,(\eta_{i}\gamma_{i})^{q_{i}}].
\]
\medskip

The proof is along similar lines to the above; as we are aiming for a slightly
weaker conclusion, we may assume from the start that each $\gamma_{j}\in\Phi$.
\medskip

\textbf{If $r$ is odd} then all the root subgroups $Y_{\omega}$ of $P$ are
one-parameter, $P$ is abelian and we set $P(2)=P$. \medskip

\textbf{If $r$ is even} then define $P(2)$ to be the product of all the
\textbf{one parameter} root subgroups $Y_{\omega}$ of $P$ together with the
$r/2$ groups
\[
B_{\omega}:= \{X_{w}(a\cdot t_{0})\ \mid\quad t_{0} \in F_{0}\},
\]
where $w$ is the root fixed by $\tau$ in an equivalence class $\omega
\subseteq\Sigma^{0}_{+}$ of type $A_{2}$ and $a$ is a fixed solution to
$a+a^{\phi}=0$.

Both $P(2)$ and $P/P(2)$ are abelian groups and modules for $\mathcal{D}\Phi$.

\medskip

We first deal with the group $P/P(2)$. It is nontrivial only if $r$ is even.
Then $P/P(2)$ is a product of its subgroups of the form
\[
\{A_{\omega}(t):=X_{v}(t)\cdot X_{u}(t^{\phi}) P(2)/P(2) \ \mid\quad t \in F
\},
\]
where the untwisted roots $v$ and $u=v^{\tau}$ span a root system
$\omega=\{u,v,u+v\}$ of type $A_{2}$ in $\Sigma^{0}_{+}$.

For $\lambda\in F_{0}^{\ast}$ let $\eta_{\lambda}$ be the inner diagonal
automorphism of $S$ induced by diag$(\lambda^{-1},\ldots,\lambda
^{-1},1,\lambda,\ldots,\lambda)$. (The unit coefficient is in the middle
position $r/2+1$ on the diagonal.) Then $(\eta_{\lambda}\gamma_{j})^{q_{j}}$
acts on each $A_{\omega}(t)\leq P/P(2)$ by
\[
t\mapsto\lambda^{\gamma_{j}+\cdots+\gamma_{j}^{q_{j}}}\cdot t^{\gamma
_{j}^{q_{j}}}.
\]
Lemma \ref{sur} part (b) gives that provided $N_{2}^{\prime}>2q(2q+1)$ then
there is a choice of $(\lambda_{1},\ldots,\lambda_{N_{2}^{\prime}})\in
F_{0}^{(N_{2}^{\prime})}$ such that the map
\[
\mathbf{t}\in F^{(N_{2}^{\prime})}\mapsto\sum_{j=1}^{N_{2}^{\prime}}\left(
\lambda_{j}^{\gamma_{j}+\cdots+\gamma_{j}^{q_{j}}}\cdot t_{j}^{\gamma
_{j}^{q_{j}}}-t_{j}\right)
\]
is surjective onto $F$. This gives (in case $r$ even)
\begin{equation}
P/P(2)=\prod_{j=1}^{N_{2}^{\prime}}[P/P(2),(\eta_{\lambda_{j}}\gamma
_{j})^{q_{j}}]. \label{twenty}%
\end{equation}

The abelian group $P(2)$ requires a little more attention since the range of
the parameter is sometimes $F$ and sometimes $F_{0}$. More precisely $P(2)$ is
a product of groups of the following three types: \medskip

Type 1: $Y_{\omega}(t_{0})=X_{w}(t_{0})$ where $\omega=\{w\}$ is a singleton
equivalence class of untwisted roots and $t_{0}$ ranges over $F_{0}$ (this
type occurs only in case $r$ is odd). \medskip

Type 2: $B_{\omega}=\{X_{w}(a\cdot t) \ \mid\ t \in F_{0}\}$, $w=v+u=w^{\tau}$
and $v,u$ span a root system $\omega$ of type $A_{2}$ (which happens only for
$r$ even). \medskip

Type 3: $Y_{\omega}(t)=X_{v}(t)X_{u}(t^{\phi}) $ where $\{u,v\}$ is an
equivalence class of untwisted roots $\omega$ of type $A_{1}\times A_{1}$ and
$t$ ranges over $F$. \medskip

Now, for $\lambda_{j} \in F^{0}$ let $\eta_{\lambda_{j}}$ be the diagonal
automorphism of $S$ defined above (with the unit coefficient omitted when $r$
is odd). Thus $\eta_{\lambda_{j}}$ acts on the parameter $t$ (or $t_{0}$)
above as multiplication by $\lambda_{j}^{2}$.

For $N_{3}^{\prime},N_{4}^{\prime}\in\mathbb{N}$ let $J_{i}$, $i=3,4$ denote
two consecutive intervals of integers of lengths $N_{i}^{\prime}$ each.

Provided $N^{\prime}_{i}$ is sufficiently big compared to $q$, then according
to Lemma \ref{sur} it is possible to find $\lambda_{j} \in F_{0}^{J_{i}}$ such
that the following two maps are surjective:%

\[
f_{3}:F_{0}^{J_{3}}\longrightarrow F_{0}\quad\mathbf{t}\mapsto\sum_{j\in
J_{1}}\left(  \lambda_{j}^{2(\gamma_{j}+\cdots+\gamma_{j}^{q_{j}})}\cdot
t_{j}^{\gamma_{j}^{q_{j}}}-t_{j}\right)
\]%
\[
f_{4}:F^{J_{4}}\longrightarrow F\quad\mathbf{t}\mapsto\sum_{j\in J_{4}}\left(
\lambda_{j}^{2(\gamma_{j}+\cdots+\gamma_{j}^{q_{j}})}\cdot t_{j}^{\gamma
_{j}^{q_{j}}}-t_{j}\right)
\]
(Note that we need part (b) of Lemma \ref{sur} for $f_{4}$.)

This implies that for $J=J_{3}\cup J_{4}$ we have
\[
P(2)=\prod_{j\in J}[P(2),(\eta_{\lambda_{j}}\gamma_{j})^{q_{j}}].
\]
Together with (\ref{twenty}) this gives the result if we take $N_{1}^{\prime
}>N_{2}^{\prime}+N_{3}^{\prime}+N_{4}^{\prime}$.

\texttt{\bigskip}

\texttt{Nikolay Nikolov}

\texttt{New College,}

\texttt{Oxford OX1 3BN,}

\texttt{UK.}

\bigskip

\texttt{Dan Segal}

\texttt{All Souls College}

\texttt{Oxford OX1 4AL}

\texttt{UK.}
\end{document}